\documentclass[11pt]{amsart}
\usepackage{amscd, amssymb, latexsym}
\usepackage{amsmath, amsthm}
\usepackage{mathrsfs}
\newtheorem{Lemma}{Lemma}[section]
\newtheorem{theorem}[Lemma]{Theorem}
\newtheorem{proposition}[Lemma]{Proposition}
\newtheorem{corollary}[Lemma]{Corollary}
\newtheorem{definition}[Lemma]{Definition}
\newtheorem{remark}[Lemma]{Remark}

\newtheorem{lemma}[Lemma]{Lemma}

\newtheorem{problem}[Lemma]{Problem}

\newcommand{\prf}[1]{\noindent {\sc Proof.} #1 \hspace*{\fill} $\Box$}

\title[Constancy of Newton polygons and isotriviality of families of curves]
{\bf Constancy of Newton polygons of $F$-isocrystals on Abelian varieties and isotriviality of families of curves}
\author[N.~TSUZUKI]{Nobuo TSUZUKI}
\address{Mathematical Institute, Tohoku University, Aza-Aoba 6-3, Aramaki, Aobaku, Sendai, 980-8578, Japan}
\email{tsuzuki@math.tohoku.ac.jp}

\subjclass[2010]{14F30 (Primary), 14H10 (Secondary)}
\keywords{Constancy of Newton polygons of $F$-isocrystals, Isotriviality of families of curves}

\begin{document}
\begin{abstract}
We prove constancy of Newton polygons of all convergent $F$-isocrystals on Abelian varieties over finite fields. 
Applying the constancy, we prove the isotriviality of projective smooth families of curves over Abelian varieties. 
We also study the isotriviality over simply connected projective smooth varieties. 
\end{abstract}

\maketitle

\section{Introduction}

The variation of Newton polygons of convergent $F$-isocrystals on algebraic varieties of characteristic $p > 0$ is mysterious. 
It may depend on the geometry of the varieties. In this paper we study the non-variation (i.e., constancy) of Newton polygons 
in projective smooth cases. 
We prove constancy of Newton polygons of all convergent $F$-isocrystals on Abelian varieties over finite fields. 
We also give several examples of projective smooth varieties such that any convergent $F$-isocrystal on the variety 
has constant Newton polygons. Applying the constancy, we prove the isotriviality of projective smooth families of such varieties. 

\subsection{Problems}
Let $k$ be a perfect field of characteristic $p$, $R$ a completely discrete valuation ring 
with residue field $k = R/\mathbf m$ where $\mathbf m$ is the maximal 
ideal of $R$, and $K$ the field of fractions of $R$ which is 
of mixed characteristic $(0, p)$. Let $\mathrm{ord}_p$ be the discrete valuation on $K$ normalized $\mathrm{ord}_p(p) = 1$. 
Let $\sigma$ be a $q$-Frobenius on $K$ for a positive power $q$ of $p$, that is, 
a continuous lift of the $q$-power Frobenius on $k$. (See the detail of Frobenius $\sigma$ in Remark \ref{frob}.) 

Let $X$ be a scheme separated of finite type over $\mathrm{Spec}\, k$, and $\mathcal M$ a convergent 
$F$-isocrystal on $X/K$ with respect to Frobenius $\sigma$ which is introduced by P.Bertehlot \cite{Be0} \cite{Be} \cite{Be1}. 
By fixing a Frobenius $\sigma$ on $K$, it is only said to be a convergent $F$-isocrystal on $X/K$ if there is no ambiguity. 
For a point $x \in X$ with a geometric point $\overline{x}$ above $x$, 
we define a polygon $\mathrm{NP}(\mathcal M, x)$, called Newton polygon of $\mathcal M$ at $x$, by 
the Newton polygon of the $F$-space $i_{\overline{x}}^\ast\mathcal M$ 
over $K(\overline{x})$ in the theory of Dieudonn\'e-Manin.  
Here $k(\overline{x})$ is a function field of $\overline{x}$ and $K(\overline{x})$ is an extension of $K$ 
with the residue field $k(\overline{x})$ as a discrete valuation field with the same valuation group, and 
$i_{\overline{x}} : \mathrm{Spec}\, k(\overline{x}) \rightarrow X$ is the canonical morphism. 
We normalize slopes by $\mathrm{ord}_p(a)/\mathrm{ord}_p(q)$ for any nonzero element $a$ of $K$ 
(note that $q$ is the $p$-power of Frobenius $\sigma$). 
We regard the application $x \mapsto \mathrm{NP}(\mathcal M, x)$ as a function 
on the scheme-theoretic points of $X$. 
The smallest slope of $\mathrm{NP}(\mathcal M, x)$ is called the initial slope of $\mathcal M$ at $x$. 

\begin{definition}\label{CNP} Let $f : X \rightarrow \mathrm{Spec}\, k$ be a morphism 
separated of finite type. 
\begin{enumerate} 
\item A convergent $F$-isocrystal $\mathcal M$ on $X/K$ is said to be constant if 
$\mathcal M \cong f^\ast\mathcal N$ for some $F$-isocrystal $\mathcal N$ on $\mathrm{Spec}\, k/K$. 
We denote  by $F\mbox{-}\mathrm{Isoc}(X/K)^{\mathrm{CST}}$ the full subcategory of constant convergent $F$-isocrystals 
in the category $F\mbox{-}\mathrm{Isoc}(X/K)$ of 
convergent $F$-isocrystals on $X/K$.
\item A convergent $F$-isocrystal $\mathcal M$ on $X/K$ is said to have constant Newton polygons 
if the application $x \mapsto \mathrm{NP}(\mathcal M, x)$ is constant on $X$. 
We denote by $F\mbox{-}\mathrm{Isoc}(X/K)^{\mathrm{CNP}}$ 
the full subcategory of convergent $F$-isocrystals with constant Newton polygons 
in $F\mbox{-}\mathrm{Isoc}(X/K)$. 
\end{enumerate}
\end{definition}

Our problems are as follow. 

\begin{problem}\label{prob} Let $X$ be a smooth and connected scheme over $\mathrm{Spec}\, k$, 
and $X_{\overline{k}} = X \times_{\mathrm{Spec}\, k} \mathrm{Spec}\, \overline{k}$ for an algebraic closure 
$\overline{k}$ of $k$. 
\begin{enumerate}
\item Is the condition $F\mbox{-}\mathrm{Isoc}(X/K) = F\mbox{-}\mathrm{Isoc}(X/K)^{\mathrm{CST}}$ 
equivalent to the triviality of geometrically etale fundamental group, i.e., $\pi_1^{\mathrm{et}}(X_{\overline{k}}) = \{ 1\}$?
Here $\pi_1^{\mathrm{et}}(-)$ denotes the etale fundamental group 
and we usually omit to indicate a base point because we do not need it in this paper. 
\item Classify varieties $X$ such that any convergent $F$-isocrystal on $X/K$ has constant Newton polygons, 
i.e., 
$$
        F\mbox{-}\mathrm{Isoc}(X/K) = F\mbox{-}\mathrm{Isoc}(X/K)^{\mathrm{CNP}}, 
$$
and study properties of such varieties. 
\end{enumerate}
\end{problem}

By definition there are natural inclusion relations 
$$
     F\mbox{-}\mathrm{Isoc}(X/K)^{\mathrm{CST}} \subset 
     F\mbox{-}\mathrm{Isoc}(X/K)^{\mathrm{CNP}} \subset 
     F\mbox{-}\mathrm{Isoc}(X/K). 
$$
When $X = \mathbb P^1_k$ is the projective line, any convergent $F$-isocystal on $X/K$ is constant. 
Indeed, there is a natural isomorphism 
$$
      \mathcal M \cong\, H^0_{\mathrm{rig}}(\mathbb P^1_k/K, \mathcal M) \otimes_K \mathcal O_{]\mathbb P^1_k[}
$$
as $F$-isocrystals. Note that the isomorphism above holds without 
assuming the existence of Frobenius structures. 
On the other hand 
there exists a projective smooth curve $C$ over a finite field 
which has a convergent $F$-isocrystal on $C/K$ 
with nonconstant Newton polygons (see the detail in Section \ref{nonisotsec}). 
Hence, even if $k$ is a finite field, 
we know there exists a projective smooth curve $C$ such that 
$$
   F\mbox{-}\mathrm{Isoc}(C/K)^{\mathrm{CST}} \subsetneq F\mbox{-}\mathrm{Isoc}(C/K)^{\mathrm{CNP}} \subsetneq 
     F\mbox{-}\mathrm{Isoc}(C/K). 
$$ 
The first inequality holds if $\pi_1^{\mathrm{et}}(C_{\overline{k}})$ 
is nontrivial by Katz-Crew equivalence \cite[Theorem 2.1]{Cr} between 
the category of 
$p$-adic continuous representations of $\pi_1^{\mathrm{et}}(C)$ and 
that of unit-root convergent $F$-isocrystals on $C/K$. 

Our main interest is the mysterious gap 
of two categories $F\mbox{-}\mathrm{Isoc}(X/K)^{\mathrm{CNP}} \subset F\mbox{-}\mathrm{Isoc}(X/K)$. 
In this paper we prove that the inclusion relations 
$$
   F\mbox{-}\mathrm{Isoc}(X/K)^{\mathrm{CST}} \subsetneq F\mbox{-}\mathrm{Isoc}(X/K)^{\mathrm{CNP}} = 
     F\mbox{-}\mathrm{Isoc}(X/K). 
$$
holds for any Abelian variety over a spectrum $\mathrm{Spec}\, k$ of a finite field $k$ (see Section \ref{conres}). 
However, the author is ignorant of classifications of projective smooth curves of genus $\geq 2$ 
in this aspect. 

Concerning to Problem \ref{prob} (1), H.Esnault and A.Shiho studied the constant problem, 
called de Jong conjecture \cite[Conjecture 2.1]{ES}, in the case 
of convergent isocrystals (without Frobenius structures). 
They proved the de Jong conjecture under certain conditions \cite{Shi} \cite{ES} \cite{ES2}. 
Moreover the constancy of geometric convergent isocrystals (see Definition \ref{gomF}) 
are proved in \cite[Theorem 1.3]{ES2}. 
As mathematical statements, the case of convergent isocrystals 
is much stronger than that of convergent $F$-isocrystals. However, our interests are the variation of Newton polygons 
and its application to geometry. 
Returning to our problems, if $X$ is proper smooth, 
then the triviality $\pi_1^{\mathrm{et}}(X_{\overline{k}}) = \{ 1\}$ is equivalent to the coincidence 
$$
     F\mbox{-}\mathrm{Isoc}(X/K)^{\mathrm{CST}} = F\mbox{-}\mathrm{Isoc}(X/K)^{\mathrm{CNP}} 
$$
of categories (Corollary \ref{van}). 
At this moment we do not know whether Problem \ref{prob} (1) is valid or not in general. 

We study basic properties of variation of Newton polygons in the section $2$, 
and give a proof of the existence of slope filtrations of 
$F$-isocrystals with constant Newton polygons in our context in Appendix \ref{esf} (see Section \ref{frobslope}). 

\subsection{Constancy results}\label{conres}

Our main result is to give an example to Problem \ref{prob} (2). 

\begin{theorem}\label{abconst0} \mbox{\rm (Theorem \ref{abconst})} 
Let $k$ be a finite field, and $X$ an Abelian variety over $\mathrm{Spec}\, k$. 
Any convergent $F$-isocrystal on $X/K$ has constant Newton polygons, i.e., 
$$
   F\mbox{-}\mathrm{Isoc}(X/K)^{\mathrm{CNP}} = F\mbox{-}\mathrm{Isoc}(X/K). 
$$
\end{theorem}

The crucial idea of the proof is as follows. Let $C$ be a projective smooth and 
geometrically connected curve of genus $\geq 1$, 
$\mathcal M$ a convergent $F$-isocrystal on $C/K$, and $D_{\mathcal M}$ 
a reduced divisor of $C$ consisting of points $x$ at which the Newton polygon $\mathrm{NP}(\mathcal M, x)$ 
is different from that at the generic point. 
Estimating the degree of $D_{\mathcal M}$ by two ways, the congruence of the $L$-function 
$L(X/k, \mathcal M; t)$ of $\mathcal M$ modulo $\mathbf m$ and the Euler-Poincar\'e formula (Proposition \ref{est}), 
we have a upper bound of the degree of $D_{\mathcal M}$ by a constant 
which depends only on the cardinal of $k$, the genus of $C$ and the 
rank of $\mathcal M$ (Theorem \ref{estjump} and Remark \ref{preest}). 

Now let $X$ be an Abelian variety over $\mathrm{Spec}\, k$. 
Let $\mathcal M$ be a convergent $F$-isocrystal on $X/K$ with non constant Newton polygons, 
and $D_{\mathcal M}$ a set of $X$ consisting of points at which 
the Newton polygon of $\mathcal M$ is different from the Newton polygon at the generic point of $X$. 
Then $D_{\mathcal M}$ is a closed subscheme in $X$ purely of codimension $1$ 
by de Jong-Oort purity theorem \cite[Theorem 4.1]{DO}. 
In this situation one can take a projective smooth and geometrically connected curve $C$ of genus $\geq 1$ in $X$ 
(we may replace $k$ by a finite extension) such that 
the degrees of the reduced divisors $(D_{[n]^\ast\mathcal M} \cap C)_{\mathrm{red}}$ 
are unbounded on $n$, where $[n]$ means the morphism of multiplication on $X$ with a positive integer $n$. 
Since the rank of $[n]^\ast\mathcal M$ is stable for $n$, 
this contradicts to the existence of bound with respect to the $k$-curve $C$.  

\subsection{Application to the isotriviality of families of curves}\label{isotsec}
We will apply our study on constancy of Newton polygons to the isotiviality problem of families of curves in characteristic $p>0$. 
At first we recall the definition and review the case of complex algebraic varieties. 

\begin{definition}\label{Isot} Let $S$ be a scheme separated of finite type over $\mathrm{Spec}\, k$. 
A smooth family $X$ over $S$ is isotrivial if, for any geometric points $s, t \in S(\overline{k})$, 
the geometric fibers $X_s$ and $X_t$ of $X$ at $s$ and $t$, 
respectively, are isomorphic to each other as schemes over $\mathrm{Spec}\, \overline{k}$. 
Here $\overline{k}$ is an algebraic closure of $k$. 
\end{definition}

In the case of complex algebraic varieties, the following results are known. 
Let $X$ be a projective smooth family of connected curves of genus $g \geq 2$ over $S$. 
If $S$ is either $\mathbb P^1$, $\mathbb C$, $\mathbb C^\times$, Abelian varieties, 
or simply connected (more generally the topological fundamental group 
$\pi_1^{\mathrm{top}}(S)$ is finite), then the family $X$ is isotrivial 
over $S$ (see the detail in \cite{DM}).

When $k$ is a field of characteristic $p>0$, 
the isotriviality of families of curves over a projective smooth curve of genus $\leq 1$ 
is known. More generally L.Szpiro proved the isotriviality of 
semistable curves over a projective smooth curve of genus $\leq 1$ \cite[Th\'eor\`eme 4]{szpiro}. 

Our result is as follows. 

\begin{theorem}\label{isotcNP0} \mbox{\rm (Theorem \ref{isotcNP})} 
Let $S$ be a projective smooth and connected scheme over $\mathrm{Spec}\, k$, 
and $f : X \rightarrow S$ a proper smooth family of connected curves. 
Suppose that any geometric convergent $F$-isocrystal $\mathcal M$ on $S/K$ (see Definition \ref{gomF}) 
has constant Newton polygons. 
Then the family $X$ over $S$ is isotrivial. 
\end{theorem}

If there exists at least a fiber which is a ordinary curve (see Definition \ref{ordi}) in the family $X$ over $S$, 
then $X$ is a family of ordinary curves over $S$ by 
the constancy hypothesis of Theorem \ref{isotcNP0}. 
In this case the family is isotrivial since the ordinary locus of coarse moduli scheme of Abelian varieties 
is quasi-affine \cite[XI, Th\'eor\`eme 5.2]{MB}. 
When there exists no ordinary point in $S$, we apply A.Tamagawa and M.Sa\"idi's works \cite{tam}, \cite{saidi} : 
there exists a finite etale covering $f' : Y' \rightarrow X'$ after a finite etale base change $S'$ of $S$ 
such that the new part $J(Y', X') = J(Y')/(f')^\ast J(X')$ of 
the relative Jacobian variety over $S'$ is ordinary. 
$J(Y, X)$ over $S$ is isotrivial and hence Torelli theorem for new parts \cite[Corollaries 4.7]{tam} 
implies the isotriviality of $X$ over $S$. 

The author is in ignorance whether the converse of the theorem above holds or not. 

\begin{corollary}\label{isotrivial} Suppose that $S$  is a projective smooth and geometrically 
connected scheme over $\mathrm{Spec}\, k$ 
satisfying one of the following:
\begin{enumerate}
\item $\pi_1^{\mathrm{et}}(S_{\overline{k}}) $ is finite. 
\item \mbox{\rm (Corollary \ref{isotab})} $S$ is an Abelian variety. 
\item $S$ is either a ruled surface over a curve of genus $\leq 1$ 
or a projective smooth surface of Kodaira dimension $0$. 
\end{enumerate}

Then any proper smooth family of connected curves over $S$ is isotrivial. 
\end{corollary}

In order to prove (1) it is sufficient to prove the assertion when $S$ is simply connected, i.e., 
the geometrically etale fundamental group $\pi_1^{\mathrm{et}}(S_{\overline{k}})$ is trivial. 
Then any geometric convergent isocrystal on $S/K$ is constant by \cite[Theorem 1.3]{ES2}, 
so that it is also constant as a convergent $F$-isocrystal. 
The isotriviality follows from Theorem \ref{isotcNP}. 
(3) follows from the classification theorem of surfaces (see \cite{Li}, for example), 
the previous (1), (2) and Propositions \ref{dominant} and \ref{birat}. 

\begin{remark} The author does not know whether 
any convergent $F$-isocrystal on a simply connected variety is constant 
or not in general. 
In the following $S$ we know that any convergent $F$-isocrystal is constant when 
the field $k$ of definition is arbitrary: 
\begin{enumerate}
\item $S$ has a projective smooth lift $\mathcal S$ over $\mathrm{Spec}\, R$ 
such that $\pi_1^{\mathrm{et}}(\mathcal S_{\overline{K}}) = \{1 \}$ 
for the geometric generic fiber $\mathcal S_{\overline{K}}$. 
\item $S$ is separably rationally connected (see the definition in \cite[IV, Definition 3.2]{Ko}). 
\end{enumerate}
Indeed, 
in the case (1) any convergent isocrystal on $S/K$ is constant, 
by rigid and complex GAGA principles of de Rham cohomologies 
and \cite{Ma} \cite{Gr} (see Introduction of \cite{Shi}), 
so that it is also constant as a convergent $F$-isocrystal. 
In the case (2) any two points are connected by rational curves \cite[IV,Theorem 3.9]{Ko} and hence 
any convergent $F$-isocrystal has constant Newton polygons. Since the simply connectedness of 
separably rationally connected varieties 
follows from the work of A.J.De Jong-J.Starr \cite{DS} and J.K\'ollar (see \cite[Corollary 3.6]{De}), the assertion follows from 
Corollary \ref{van}. 
\end{remark}

\vspace*{3mm}

\noindent
{\bf Acknowledgments.} The author thanks Professor Takao Yamasaki and Professor Jeng-Daw Yu 
for useful discussions. He also thanks Professor Yifan Yang who told the author 
important examples of Shimura curves. 
The examples inspired the author to study the constancy problem. 
The author is supported by the Grant-in-Aid for Exploratory Research, 
Japan Society for the promotion of Science.

\section{Newton polygon}

In this section we study several properties of Newton polygons. 

\subsection{Frobenius $\sigma$}

Let $k$ be a perfect field of characteristic $p > 0$, $R$ a completely discrete valuation ring 
with residue field $k = R/\mathbf m$, and $K$ the field of fractions of $R$ which is 
of mixed characteristic $(0, p)$. 
We fix a discrete valuation $\mathrm{ord}_p$ on $K$ (and its extension as valuation fields) which is normalized 
by $\mathrm{ord}_p(p) = 1$. 
Let $\sigma$ be a $q$-Frobenius on $K$ for a positive power $q$ of $p$, that is, 
a continuous lift of the $q$-power Frobenius on $k$. 

When we consider Katz-Crew equivalence between continuous $p$-adic representations 
and unit-root $F$-isocrystals \cite[Theorem 2.1]{Cr}, we assume that 
\begin{enumerate}
\item[(i)] $\mathbb F_q \subset k$ and 
\item[(ii)] $K_\sigma \otimes_{W(\mathbb F_q)}W(k) \cong K$, 
\end{enumerate}
where 
$\mathbb F_q$ is the finite field of $q$-elements, 
$K_\sigma$ is the $\sigma$-invariant 
subfield of $K$ and $W(k)$ is the ring of Witt-vectors with coefficients in $k$. 

\begin{remark}\label{frob} 
The above hypotheses (i) and (ii) on $K$ and $\sigma$ 
always hold if we replace $K$ 
by a finite unramified extension $K'$ of $K$. Indeed, for any $q$-Frobenius $\sigma$, there 
is a uniformizer $\pi$ of the $p$-adic completion $\widehat{K}^{\mathrm{ur}} = K \otimes_{W(k)}W(\overline{k})$ 
of a maximal unramified extension $K^{\mathrm{ur}}$ of $K$ such that 
which is $\sigma$-fixed, i.e., $\sigma(\pi) = \pi$. 
Here $\sigma$ also denotes the unique extension of Frobenius $\sigma$ on $K^{\mathrm{ur}}$. 
Since $(\widehat{K}^{\mathrm{ur}})_\sigma$ is a finite extension of $\mathbb Q_p$, $\pi$ is algebraic over $K$. 
If we put $K'= K(\widehat{K}^{\mathrm{ur}})_\sigma(W(\mathbb F_q)[1/p])$ which is a composite of fields, 
then $K'$ is unramified over $K$ by $K' \subset \widehat{K}^{\mathrm{ur}}$ 
and $K_\sigma \otimes_{W(\mathbb F_q)}W(k') \cong K'$ for some finite extension $k'$ of $k$. 
\end{remark}

\subsection{Variation of Newton polygons}

Let $M$ be an $F$-space over $K$ with respect to Frobenius $\sigma$, that is, 
a $K$-vector space $M$ of finite dimension with a $\sigma$-linear bijection 
$$
        F : M \rightarrow M 
$$
which is called Frobenius. 
Suppose that the residue field $k$ is algebraically closed and fixed a uniformizer 
$\pi$ of $K_\sigma$. The category of $F$-spaces is semi-simple with simple objects 
$$
       E_{s/r} = K[F]/(F^r - \pi^s)\, \hspace*{5mm} \mbox{for}\, \, 
       (r, s) \in \mathbb Z_{> 0} \times \mathbb Z\,\,\mbox{with}\, (r, s) = 1
$$
where the relation $Fa = \sigma(a)F$ holds for any $a \in K$ in $K[F]$, by Dieudonn\'e-Manin's classification. 
The rational number $s\mathrm{ord}_p(\pi/q)/r$ is called slope of $E_{r, s}$ and its rank over $K$ is $r$. 
If $M = \oplus_iE_{r_i/s_i}^{m_i}$ with $s_1/r_1<s_2/r_2<\cdots<s_l/r_l$, then the Newton polygon of $M$ 
is the lower convex full of points 
$$
   \begin{array}{l}
     (0, 0), (m_1r_1, m_1s_1\mathrm{ord}_p(\pi/q)), 
     (m_1r_1+m_2r_2, (m_1s_1+m_2s_2)\mathrm{ord}_p(\pi/q)), \cdots, \\
    \hspace*{45mm}  (m_1r_1+\cdots +m_lr_l, (m_1s_1+\cdots+m_ls_l)\mathrm{ord}_p(\pi/q)).
     \end{array}
$$
Note that if $F = \pi\sigma$ of rank $1$, then it is of slope $\mathrm{ord}_p(\pi/q)$. Our definition of slopes 
is stable under the extension of scalar field $K$ and the change of power of Frobenius $\sigma$. 

Let $X$ be a scheme separated of finite type over $\mathrm{Spec}\, k$ and $\mathcal M$ a convergent $F$-isocrystal 
on $X/K$. For a point $x \in X$ (not necessary a closed point) and for a geometric point $\overline{x}$ with a natural morphism 
$i_{\overline{x}} : \overline{x} \rightarrow X$, we define the Newton polygon of $\mathcal M$ at $x$ by 
that of $F$-space $i_{\overline{x}}^\ast\mathcal M$ 
over $K(\overline{x})$. Here $k(\overline{x})$ is the field of functions of $\overline{x}$, 
and $K(\overline{x})$ is the extension of $K$ as complete discrete valuation fields 
such that the residue field is $k(\overline{x})$. It is independent of the choice of geometric point $\overline{x}$ above $x$. 
We denote Newton polygon of $\mathcal M$ at $x$ by $\mathrm{NP}(\mathcal M, x)$. Then 
$$
       x \in X\,\, \mapsto\,\, \mathrm{NP}(\mathcal M, x)
$$ 
is a function on $X$. 

Let $\mathrm{NP}_1$ and $\mathrm{NP}_2$ be Newton polygons. We say that $\mathrm{NP}_1$ is above $\mathrm{NP}_2$, 
denote it by $\mathrm{NP}_1 \prec \mathrm{NP}_2$, 
if $\mathrm{NP}_1$ and $\mathrm{NP}_2$ has same endpoints and all polygons of $\mathrm{NP}_1$ is 
upper than or equal to that of $\mathrm{NP}_2$. 

\begin{proposition}\label{sp} \mbox{\rm (Grothendieck's specialization theorem \cite[Theorem 2.3.1]{Kat})} 
With the notation as above we have the following. 
\begin{enumerate}
\item Let $x, y \in X$. If $x$ is a specialization of $y$, then 
$\mathrm{NP}(\mathcal M, x) \prec \mathrm{NP}(\mathcal M, y)$. 
\item Suppose $X$ is irreducible. 
If $\gamma$ is the initial slope of $\mathcal M$ at the generic point of $X$, then the set 
$$
     U = \{ x \in X\, |\, \mbox{\rm the initial slope of $\mathcal M$ at $x$ is $\eta$.} \}
$$
is open in $X$. 
\end{enumerate}
\end{proposition}

Another important property of variation of Newton polygon 
is the purity theorem of Newton polygons by A.de Jong and F.Oort. 

\begin{theorem}\label{Jort} \mbox{\rm (\cite[Theorem 4.1]{DO})} 
Let $X$ be a smooth irreducible scheme separated of finite type and 
$\eta$ a generic point of $X$. For a convergent $F$-isocrystal $\mathcal M$ 
on $X/K$, any generic point of the set-theoretical 
complement of the open subscheme 
$$
  U_{\mathcal M} = \{ x \in X\, |\, \mathrm{NP}(\mathcal M, x) = \mathrm{NP}(\mathcal M, \eta) \}
$$
is of codimension $1$ in $X$. 
\end{theorem}

\subsection{Slope filtrations}\label{frobslope}

\begin{definition}\label{sfdef} Let $X$ be a scheme separated of finite type over $\mathrm{Spec}\, k$, and 
$\mathcal M$ a convergent $F$-isocrystal on $X/K$. 
\begin{enumerate}
\item $\mathcal M$ is said to be isoclinic of slope $\gamma \in\mathbb Q$ 
if, for any geometric point $\overline{x}$ of $X$ with a natural morphism $i_{\overline{x}} : \overline{x} \rightarrow X$, 
the $F$-space $i_{\overline{x}}^\ast\mathcal M$ over $K(\overline{x})$ is a direct sum of copies of 
$E_{\gamma/\mathrm{ord}_p(\pi/q)}$. 
\item An increasing filtration $\{S_\lambda \mathcal M\}_{\lambda \in \mathbb Q}$ of $\mathcal M$ as convergent 
$F$-isocrystals on $X/K$ is called the slope filtration if it satisfies the following conditions
\begin{enumerate}
\item[(i)] $S_\lambda \mathcal M = 0$ for $\lambda <\hspace*{-1.5mm}< 0$, 
$S_\lambda \mathcal M = \mathcal M$ for $\lambda >\hspace*{-1.5mm}> 0$ and 
$S_{\lambda+} \mathcal M = S_\lambda \mathcal M$ for any $\lambda$; 
\item[(ii)] $S_\lambda \mathcal M/S_{\lambda -} \mathcal M$ is 
either $0$ or isoclinic of slope $\lambda$ for any $\lambda$, 
\end{enumerate}
where $S_{\lambda -} \mathcal M = \cup_{\mu < \lambda} S_\mu \mathcal M$ 
and $S_{\lambda +} \mathcal M = \cap_{\mu > \lambda} S_\mu \mathcal M$. 
\end{enumerate}
\end{definition}

The existence of slope filtrations of $F$-crystals 
on $\mathrm{Spec}\, k[\hspace*{-0.4mm}[t]\hspace*{-0.4mm}]$ having constants Newton polygons up to isogenies 
was proved in \cite[Corollary 2.6.3]{Kat}, and 
the case of unipotent convergent $F$-isocrystals 
was proved in \cite[Th\'eor\`eme  3.2.3]{CLpente}. 
The existence of generic slope filtration for convergent $F$-isocrystals is proved in \cite[Proposition 5.8]{Shi0}. 
In the general case the existence of slope filtration is stated in \cite[Corollary 4.2]{Ke2}. 
The author does not find the detail of proof in the literature. So 
we prove the following theorem in Appendix \ref{esf}. 

\begin{theorem}\label{fil} 
Let $X$ be a smooth scheme separated of finite type over $k$, and $\mathcal M$ 
a convergent $F$-isocrystal on $X/K$. Suppose that the initial slope (i.e., the smallest slope) of $\mathcal M$ 
of the generic point of $X$ is $\gamma$ and the multiplicity of slope $\gamma$ 
of the induced $F$-isocrystal $i_x^\ast\mathcal M, (x \in X)$ is constant on $X$, say the rank is $r$. 
Then there exists a convergent sub $F$-isocrystal $\mathcal L$ of $\mathcal M$ on $X/K$ 
which is isoclinic of slope $\gamma$ and of rank $r$. 
\end{theorem}

\begin{corollary}\label{sf} Let $X$ be a smooth scheme separated of finite type over $\mathrm{Spec}\, k$, 
and $\mathcal M$ a convergent $F$-isocrystal on $X/K$ which has constant Newton polygons. 
Then $\mathcal M$ admits a unique slope filtration $\{S_\lambda \mathcal M\}_{\lambda \in \mathbb Q}$. 
\end{corollary}

\subsection{Properties of constancy of Newton polygons}

Let $X$ be a scheme separated of finite type over $\mathrm{Spec}\, k$, 
and $\mathcal M$ a convergent $F$-isocrystal on $X/K$. 
We define the convergent cohomology of $X$ by 
$$
      H^i_{\mathrm{conv}}(X/K, \mathcal M) = I\hspace*{-1mm}H^i(]X[_{\mathcal P}, 
      \mathcal M \otimes \Omega_{]X[_{\mathcal P}/K}^\bullet).
$$
Here we take a closed immersion $X \rightarrow \mathcal P$ into 
a smooth formal scheme $\mathcal P$ over $\mathrm{Spf}\, R$,  
$]X[_{\mathcal P}$ is the associated rigid analytic tube \cite[Section 1.1]{Be}, 
and $I\hspace*{-1mm}H^i(]X[_{\mathcal P}, \mathcal M \otimes \Omega_{]X[_{\mathcal P}/K}^\bullet)$ is 
the hypercohomology of the de Rham complex $\mathcal M \otimes \Omega_{]X[_{\mathcal P}/K}^\bullet$ 
associated to the convergent isocrystal $\mathcal M$. 
The convergent cohomology is independent of the choice of the closed immersion $X \rightarrow \mathcal P$. 
In \cite{CT0} \cite{tsu} 
the convergent cohomology 
$H^i_{\mathrm{conv}}(X/K, \mathcal M)$ is denoted by $H^i_{\mathrm{rig}}((X, X)/K, \mathcal M)$, 
the rigid cohomology of $X$ overconvergent along $\emptyset$, namely, we do not consider the overconvergent 
regularity along boundary. 
When $\mathcal M = \mathcal O_{]X[}$ (the unit $F$-isocrystal), we simply denote 
the convergent cohomology by $H^i_{\mathrm{conv}}(X/K)$. 
If $X$ is proper over $\mathrm{Spec}\, k$, then the convergent cohomology $H^i_{\mathrm{conv}}(X/K, \mathcal M)$ 
is nothing but the rigid cohomology $H^i_{\mathrm{rig}}(X/K, \mathcal M)$. 
The convergent isocrystal is furnished with a $\sigma$-linear homomorphism 
$$
        F : H^i_{\mathrm{conv}}(X/K, \mathcal M) \rightarrow H^i_{\mathrm{conv}}(X/K, \mathcal M)
$$
In general $H^i_{\mathrm{conv}}(X/K, \mathcal M)$ 
is not of finite dimension over $K$ and $F$ does not acts on $V$ bijectively except $i = 0$. 
For a $K$-vector space with $\sigma$-linear homomorphism  
$F : V \rightarrow V$, we define a $K$-space $V_{\mathrm{fin}}$ by the subspace of 
$V$ consisting of elements $w$ 
such that there exist elements $a_1, a_2, \cdots, a_n \in K$ with $a_n \ne 0$ 
satisfying $F^nw + a_1F^{n-1}w + \cdots + a_nw = 0$. In the case of rigid cohomology 
we have an equality $H^i_{\mathrm{rig}}(X/K, \mathcal M)_{\mathrm{fin}} = 
H^i_{\mathrm{rig}}(X/K, \mathcal M)$. 

\begin{proposition}\label{Abel} Let $0 \rightarrow \mathcal L \rightarrow \mathcal M \rightarrow \mathcal N \rightarrow 0$ 
be an exact sequence of $F\mbox{-}\mathrm{Isoc}(X/K)$. 
\begin{enumerate}
\item $\mathcal M$ is an object in $F\mbox{-}\mathrm{Isoc}(X/K)^{\mathrm{CNP}}$ 
if and only if both $\mathcal L$ and $\mathcal N$ are so.
\item If $\mathcal M$ is an object in $F\mbox{-}\mathrm{Isoc}(X/K)^{\mathrm{CST}}$, 
then so are both $\mathcal L$ and $\mathcal N$. Suppose furthermore that 
$H_{\mathrm{conv}}^1(X/K)_{\mathrm{fin}}  = 0$. Then the converse holds.
\end{enumerate}
In particular, the category 
$F\mbox{-}\mathrm{Isoc}(X/K)^{\mathrm{CNP}}$ (resp. 
$F\mbox{-}\mathrm{Isoc}(X/K)^{\mathrm{CST}}$) is an Abelian subcategory 
of $F\mbox{-}\mathrm{Isoc}(X/K)$. 
\end{proposition}

\vspace*{3mm}

We give several properties of constancy of Newton polygons. 

\begin{proposition}\label{dominant} 
Let $f : X \rightarrow Y$ be a morphism of schemes separated of finite type over $\mathrm{Spec}\, k$. 
\begin{enumerate}
\item Suppose that $Y$ is smooth and the set-theoretical complement of the image $f(X)$ in $Y$ 
is of codimension $\geq 2$. 
If any convergent $F$-isocrystal on $X/K$ has constant Newton polygons, then 
the same holds for any convergent $F$-isocrystal on $Y/K$. 
\item Suppose that the morphism $f$ is finite etale. 
If any convergent $F$-isocrystal on $Y/K$ has constant Newton polygons, then 
the same holds for any convergent $F$-isocrystal on $X/K$. 
\end{enumerate}
\end{proposition}

\prf{(1) follows from Theorem \ref{Jort}.

(2) Let $f_{\mathrm{rig} \ast} : F\mbox{-}\mathrm{Isoc}(X/K) \rightarrow F\mbox{-}\mathrm{Isoc}(Y/K)$ 
be the push forward induced by the finite etale morphism $f : X \rightarrow Y$. 
For a convergent $F$-isocrystal $\mathcal M$ on $X/K$ and a point $x \in X$, the Newton polygon of 
 $f_{\mathrm{rig} \ast}\mathcal M$ at a point $f(x)$ is a $\mathrm{deg}(f)$ time Newton polygon $\mathrm{NP}(\mathcal M, x)$. 
}

\begin{proposition}\label{birat} 
Let $X$ and $Y$ be projective smooth and connected 
schemes over $\mathrm{Spec}\, k$. Suppose that $X$ and $Y$ 
are binational over $\mathrm{Spec}\, k$. 
Any convergent $F$-isocrystal on $X/K$ 
is constant (resp. has constant Newton polygons) if and only 
if the same holds for any convergent $F$-isocrystal on $Y/K$. 
\end{proposition}

\prf{Since the set of fundamental points of binational transformations is closed of codimension $\geq 2$ 
\cite[V, Lemma 5.1]{Ha}, 
the assertion follows from K.Kedlaya's extension theorem \cite[Proposition 5.3.3]{Ke} (via 
the extension of overconvergent $F$-isocrystals) and Theorem \ref{Jort}. 
}

\subsection{Geometric $F$-isocrystals}

\begin{definition}\label{gomF} Let $S$ be a scheme separated of finite type over $\mathrm{Spec}\, k$. 
A convergent isocrystal (resp. $F$-isocrystal) $\mathcal M$ on $S/K$ is geometric 
if there exists a proper smooth morphism $f : X \rightarrow S$ 
such that $\mathcal M$ is isomorphic to a subquotient of the relative rigid cohomology 
$R^if_{\mathrm{rig} \ast}\mathcal O_{]X[}$ of $X$ over $S$ for some integer $i$ 
as a convergent isocrystal (resp. $F$-isocrystal). 
\end{definition}

Since $f : X \rightarrow S$ is proper smooth, the rigid cohomology is nothing but the convergent cohomology 
of A.Ogus in \cite{Og} and the finiteness theorem below holds 
(see the detail around Berthelot's conjecture of coherency of relative rigid cohomology in \cite{La}). 
Hence the above definition makes sense. 

\begin{theorem}\label{coh} \mbox{\rm (\cite[Sect. 3]{Og})} With the notation as above, 
suppose that $S$ is smooth over $\mathrm{Spec}\, k$. 
Then for any convergent $F$-isocrystal $\mathcal M$ on $X/K$, 
the relative rigid cohomology $R^if_{\mathrm{rig} \ast}\mathcal M$ is a convergent $F$-isocrystal on $S/K$. 
\end{theorem}

We will use the following proposition to reduce the problem in the case where 
the field $k$ of definition is a finite field. This similar argument is discussed in \cite[Section 5.4]{ES2}. 

\begin{proposition}\label{geom} Let $\kappa$ be a perfect field, 
$\mathcal T$ a smooth integral scheme separated of finite type over 
$\mathrm{Spec}\, \kappa$, and $\mathfrak f : \mathcal X \rightarrow \mathcal S$ 
a proper smooth morphism of smooth schemes separated of finite type 
over $\mathcal T$. Let $k$ be a perfect field which includes 
the function field $\kappa(\mathcal T)$ of $\mathcal T$, and 
$f : X \rightarrow S$ the base change of $\mathfrak f : \mathcal X \rightarrow \mathcal S$ 
by the natural morphism $\mathrm{Spec}\, k \rightarrow \mathcal T$. 
Suppose that, for any closed point $t$ of $\mathcal T$ over $\mathrm{Spec}\, \kappa$, 
the relative rigid cohomology $R^i\mathfrak f_{t, \mathrm{rig}\ast}\mathcal O_{]{\mathcal X}_t[}$ 
on $\mathcal S_t/K(t)$ has constant Newton polygons, where 
$\kappa(t)$ is the function field of $t$ and 
$K(t) = K_\sigma\otimes_{W(\mathbb F_q)}W(\kappa(t))$. Then the 
relative rigid cohomology $R^if_{\mathrm{rig}\ast}\mathcal O_{]X[}$ 
on $S/K$ has constant Newton polygons. 
\end{proposition}

\prf{We may assume that $S$ is irreducible. 
Let $\mathcal U$ be an open subscheme of $\mathcal S$ consisting of points 
at which the Newton polygon of the relative 
rigid cohomology $R^i\mathfrak f_{\mathrm{rig}\ast}\mathcal O_{]{\mathcal X}[}$ 
coincides with that at the generic point of $\mathcal S$, and hence with that at the generic point 
of $S$. The open subscheme $\mathcal U$ 
is defined over $\mathrm{Spec}\, \kappa$ by Grothendieck's 
specialization theorem. If $\mathfrak g : \mathcal S \rightarrow \mathcal T$ 
denotes the structure morphism, then $\mathfrak g(\mathcal U)$ 
is an open subscheme of $\mathcal T$ over $\mathrm{Spec}\, \kappa$ and 
$\mathfrak g^{-1}(\mathfrak g(\mathcal U)) = \mathcal U$. 
Hence $S \subset \mathcal U \times_{\mathcal T}\mathrm{Spec}\, k$, and 
$R^if_{\mathrm{rig}\ast}\mathcal O_{]X[}$ has constant Newton polygons on $S$. 
Note that, in our case the arbitrary base change theorem 
holds since the relative cohomology is coherent \cite[Corollary 2.3.3]{tsu}. 
}

\subsection{Constant and constant Newton polygons}

\begin{proposition}\label{equiv}
Let $X$ be a smooth scheme separated of finite type over $\mathrm{Spec}\, k$. Then the following conditions (i) and (ii) 
are equivalent.
\begin{enumerate}
\item[(i)] $F\mbox{-}\mathrm{Isoc}(X/K)^{\mathrm{CST}} = F\mbox{-}\mathrm{Isoc}(X/K)^{\mathrm{CNP}}$. 
\item[(ii)] $\pi_1^{\mathrm{et}}(X_{\overline{k}}) = \{ 1 \}$ and 
$H^1_{\mathrm{conv}}(X/K)_{\mathrm{fin}} = 0$ (see the definition before Proposition \ref{Abel}). 
\end{enumerate}
\end{proposition}

\prf{(i) $\Rightarrow$ (ii) follows from Proposition \ref{Abel} and Lemma \ref{ext} below. 

Now we prove (ii) $\Rightarrow$ (i). There exists a natural commutative diagram 
$$
     \begin{array}{ccc}
     \mathrm{Rep}_{K_\sigma}(\mathrm{Gal}(\overline{k}/k)) &\overset{\cong}{\longrightarrow} 
&F\mbox{-}\mathrm{Isoc}(\mathrm{Spec}\, k/K)^0 \\
     \downarrow & &\downarrow \\
     \mathrm{Rep}_{K_\sigma}(\pi_1^{\mathrm{et}}(X)) &\overset{\cong}{\longrightarrow} &F\mbox{-}\mathrm{Isoc}(X/K)^0. \\
\end{array}
$$ 
by Katz-Crew's equivalence between the category of $p$-adic continuous 
representations of etale fundamental groups 
and the category of unit-root (i.e., isoclinic of slope $0$) convergent 
$F$-isocrystals \cite[Theorem 2.1]{Cr}. Since there exists an exact sequence 
$$
      1 \rightarrow \pi_1^{\mathrm{et}}(X_{\overline{k}}) \rightarrow \pi_1^{\mathrm{et}}(X) 
      \rightarrow \mathrm{Gal}(\overline{k}/k) \rightarrow 1 
$$
of etale fundamental groups, the triviality  $\pi_1^{\mathrm{et}}(X_{\overline{k}}) = \{ 1\}$ implies 
the natural inverse image functor 
$$
F\mbox{-}\mathrm{Isoc}(\mathrm{Spec}\, k/K)^0 
\rightarrow F\mbox{-}\mathrm{Isoc}(X/K)^0
$$
is an equivalence. Let $\mathcal M$ be a convergent $F$-isocrystal on $X/K$ 
which is isoclinic, and take a finite extension $K'$ of $K$ 
whose valuation group contains $\gamma$. Let $k'$ be a residue field of $K'$, 
$X' = X \times_{\mathrm{Spec}\, k} \mathrm{Spec}\, k'$, and 
$\mathcal M'$ the inverse image of $\mathcal M$ on $X'/K'$. Then 
$\mathcal M'$ is a tensor product of a rank $1$ constant object and a unit-root object. 
Then we have an isomorphism 
$$
     H_{\mathrm{conv}}^0(X/K, \mathcal M) \otimes_K K' 
     \cong H_{\mathrm{conv}}^0(X'/K', \mathcal M')
$$
of $K'$-spaces of dimension $\mathrm{rank}\, \mathcal M$ by the unit-root case, 
and hence $\mathcal M$ is constant. 

Let $\mathcal M$ be an object in $F\mbox{-}\mathrm{Isoc}(X/K)^{\mathrm{CNP}}$ 
with a slope filtration $\{ S_\lambda\mathcal M\}$ (Corollary \ref{sf}). We have only to prove 
the slope filtration is split. This follows from the exact sequence
$$
     0 \rightarrow \mathcal L \rightarrow \mathcal M \rightarrow \mathcal N \rightarrow 0
$$
is split when $\mathcal L$ and $\mathcal N$ are isoclinic. 
Indeed, since $\mathcal L = \mathcal O_{]X[}\otimes_KH^0_{\mathrm{conv}}(X/K, \mathcal L)$, 
we have 
$$
       H^1_{\mathrm{conv}}(X/K, \mathcal L) \cong 
       H^1_{\mathrm{conv}}(X/K) \otimes_KH^0_{\mathrm{conv}}(X/K, \mathcal L). 
$$
Hence the vanishing $H^1_{\mathrm{conv}}(X/K)_{\mathrm{fin}} = 0$ implies the splitting.}

\vspace*{3mm}

The cohomological interpolation of Hom and Ext for convergent isocrystals 
in \cite[Proposition 2.2.7]{Be} and \cite[Proposition 1.2.2]{CL} implies the lemma below. 

\begin{lemma}\label{ext} Let $\mathcal M, \mathcal N$ be convergent 
$F$-isocrystals on $X/K$, and let us denote a $K_\sigma$-space of 
homomorphisms as convergent 
$F$-isocrystals (resp. a $K_\sigma$-space of extension classes of convergent $F$-isocrystals) 
on $X/K$ by 
$\mathrm{Hom}_{\tiny \mbox{\rm $F$-Isoc}}(\mathcal M, \mathcal N)$ 
(resp. $\mathrm{Ext}^1_{\tiny \mbox{\rm $F$-Isoc}}(\mathcal M, \mathcal N)$). Then there is 
an exact sequence of $K_\sigma$-spaces:
$$
  \begin{array}{ccccccc}
  0 &\rightarrow &\mathrm{Hom}_{\tiny \mbox{\rm $F$-Isoc}}(\mathcal M, \mathcal N) &\rightarrow 
  &H^0_{\mathrm{conv}}(X/K, \mathcal M^\vee \otimes\mathcal N) 
  &\overset{1-F}{\rightarrow} &H^0_{\mathrm{conv}}(X/K, \mathcal M^\vee \otimes\mathcal N) \\
   &\rightarrow &\mathrm{Ext}^1_{\tiny \mbox{\rm $F$-Isoc}}(\mathcal M, \mathcal N) &\rightarrow 
  &H^1_{\mathrm{conv}}(X/K, \mathcal M^\vee \otimes\mathcal N) &\overset{1-F}{\rightarrow} 
  &H^1_{\mathrm{conv}}(X/K, \mathcal M^\vee \otimes\mathcal N). 
  \end{array}
$$
Here $\mathcal M^\vee$ is the dual of $\mathcal M$. 
If furthermore the residue field $k$ is algebraically closed, then 
the homomorphism $H^0_{\mathrm{conv}}(X/K, \mathcal M^\vee \otimes\mathcal N) 
\overset{1-F}{\rightarrow} H^0_{\mathrm{conv}}(X/K, \mathcal M^\vee \otimes\mathcal N)$ 
is surjective. 
\end{lemma}

\vspace*{3mm}

The author does not know whether $\pi_1^{\mathrm{et}}(X_{\overline{k}}) = \{ 1 \}$
implies $H^1_{\mathrm{conv}}(X/K)_{\mathrm{fin}} = 0$ or not in general because 
the convergent cohomology is huge in general. 
Esnault and Shiho proved the following theorem by a comparison with $\ell$-adic cohomology theory. 
Hence $H^1_{\mathrm{conv}}(X/K)_{\mathrm{fin}} = 0$ at least $X$ is proper. 

\begin{theorem}\label{vann} \mbox{\rm (\cite[Theorem 5.1]{ES2})}  
Let $X$ be a smooth connected scheme over $\mathrm{Spec}\, k$. 
Suppose that either $X$ is proper over $\mathrm{Spec}\, k$ or $p \geq 3$.  
Then the triviality $\pi_1^{\mathrm{et, ab}}(X_{\overline{k}}) = \{ 1 \}$ implies the vanishing 
$H^1_{\mathrm{rig}}(X/K) = 0$. Here $\pi_1^{\mathrm{et, ab}}(X_{\overline{k}})$ is the maximal Abelian quotient 
of geometric etale fundamental group $\pi_1^{\mathrm{et}}(X_{\overline{k}})$. 
\end{theorem}

\begin{corollary}\label{van} 
Let $X$ be a smooth connected scheme over $\mathrm{Spec}\, k$. 
Suppose that $X$ is proper over $\mathrm{Spec}\, k$. 
The equivalence $F\mbox{-}\mathrm{Isoc}(X/K)^{\mathrm{CST}} = F\mbox{-}\mathrm{Isoc}(X/K)^{\mathrm{CNP}}$ 
of categories holds if and only if $\pi_1^{\mathrm{et}}(X_{\overline{k}}) = \{ 1 \}$. 
\end{corollary}

\section{Constancy of Newton polygons on Abelian varieties} 
In this section we prove any $F$-isocrystal on an Abelian variety over a finite field 
has constant Newton polygons. 

\subsection{An estimate of number of jumping points}

Suppose $k$ is a finite field of $q$ elements in this section. 
Let us denote an algebraic closure of $k$ (resp. the $p$-adic completion of a maximal unramified extension $R^{\mathrm{ur}}$ of $R$, 
resp. the field of fractions of $\widehat{R}^{\mathrm{ur}}$) 
by $\overline{k}$ (resp. $\widehat{R}^{\mathrm{ur}}$, resp. $\widehat{K}^{\mathrm{ur}}$). 
Let $\sigma$ be a $q$-Frobenius on $K$ such that $(\widehat{K}^{\mathrm{ur}})_\sigma = K$. 
When we replace $k$ by a finite extension $k'$ of $q^r$ elements, then we also change 
the Frobenius $\sigma$ by $\sigma^r$. 

Let $C$ be a projective smooth and geometrically connected curve of genus $g$ over $\mathrm{Spec}\, k$ 
with $p$-rank $e$, i.e., 
$$
   e = \mathrm{rank}_{\mathbb F_p}J(C)[p](\overline{k}) = \mathrm{dim}_KH^1_{\mathrm{rig}}(C/K)^0.
$$
Here $J(C)$ is the Jacobian variety of $C$, $J(C)[p]$ is the subgroup scheme of $J(C)$ 
which is the kernel of the multiplication by $p$, $H^i_{\mathrm{rig}}(C/K)$ 
is the $i$-th rigid cohomology of $C/K$, and $H^i_{\mathrm{rig}}(C/K)^0$ 
is the unit-root subspace of the $F$-space $H^1_{\mathrm{rig}}(C/K)$. 
Then $0 \leq e \leq g$. 

\begin{proposition}\label{est} With the notation as above, suppose $g \geq 1$. 
Let $\mathcal M$ be a convergent $F$-isocrystal on $C/K$ 
which satisfies the following conditions:
\begin{enumerate}
\item[(i)] The initial slope of $\mathcal M$ at the generic point of $C$ is $0$ of multiplicity $1$; 
\item[(ii)] Let $U = \{ x \in C\, |\,  \mbox{the initial slope of $\mathcal M$ at x is 0}\}$ be an open subscheme of $C$ 
with the open immersion $j_U : U \rightarrow X$, 
$\mathcal L$ the rank $1$ unit-root convergent sub-$F$-isocrystal 
of $j^\ast\mathcal M$ on $U/K$ by Theorem \ref{fil}, and 
$\rho : \pi_1^{\mathrm{et}}(U) \rightarrow R^\times$ 
the $p$-adic representation corresponding to $\mathcal L$ by Katz-Crew's equivalence. Then 
the representation $\rho$ satisfies 
$$
    \rho\, \equiv\, 1\, \, (\mathrm{mod}\, {\mathbf m}R). 
$$
\end{enumerate}
If $Z = C \setminus U$ is a reduced divisor of $C$ over $\mathrm{Spec}\, k$, then we have an inequality 
$$
       e + \mathrm{deg}\, Z \leq 1 + 2(g-1)\mathrm{rank}\mathcal M. 
$$
\end{proposition} 

\prf{If $U = C$, then 
there is nothing to prove. So we suppose $U \ne C$ and hence the rank of $\mathcal M$ is greater than or equal to $2$. 
Let us calculate the $L$-function 
$$
    L(C/k, \mathcal M; t) = 
     \displaystyle{\prod_{x : \mathrm{closed\, \, points\,\, of}\, C}
     \mathrm{det}(1-F_x^{\mathrm{deg}(x)}t^{\mathrm{deg}(x)}; i_x^\ast\mathcal M)^{-1}}
$$
of $\mathcal M$ modulo $\mathbf m$ by two ways. 
Here $\mathrm{deg}(x)$ is the degree of the function field $k(x)$ of $x$ over $k$ 
and $i_x : x \rightarrow C$ be the canonical morphism. 
At a closed point $x$ of $U$ all the $p$-adic valuations of Frobenius eigenvalues of $i_x^\ast\mathcal M$ 
except one eigenvalue which is $1$ modulo $\mathbf m$ are positive. 
At a closed point $x$ in $Z$ all the $p$-adic valuations of Frobenius eigenvalues of $i_x^\ast\mathcal M$ 
are positive. Hence we have 
$$
  \begin{array}{lll}
     L(C/k, \mathcal M; t)\, &\equiv\, 
     &\displaystyle{\prod_{x : \mathrm{closed\, \, points\,\, of}\, U} (1- t^{\mathrm{deg}(x)})^{-1}} \\
     &\equiv\, &\mathrm{Zeta}(C/k; t)\mathrm{Zeta}(Z/k; t)^{-1} \\
     &\equiv\, 
     &\displaystyle{\frac{\mathrm{det}(1-Ft; H^1_{\mathrm{rig}}(C/K))\mathrm{Zeta}(Z/k; t)^{-1}}{1-t}}\, \, (\mathrm{mod}\, \, \mathbf m)
     \end{array}
     \leqno{(\ast)}
$$
in $k[\hspace*{-0.4mm}[t]\hspace*{-0.4mm}]$, 
where $\mathrm{Zeta}(C/k; t)$ (resp. $\mathrm{Zeta}(Z/k; t)$) is the zeta function of $C$ (resp. $Z$) over $k$. 
Since 
$$
     \begin{array}{l}
       \mathrm{deg}(\mathrm{det}(1-Ft; H^1_{\mathrm{rig}}(C/K))\, (\mathrm{mod}\, \mathbf m)) = e \\
       \mathrm{deg}(\mathrm{Zeta}(Z/k; t)^{-1}\, (\mathrm{mod}\, \mathbf m)) = \mathrm{deg}(Z)
     \end{array}
$$
in $k[t]$, we have 
$$
  \mathrm{deg}((1-t)L(C/k, \mathcal M; t)\,  (\mathrm{mod}\, \, \mathbf m)) = e + \mathrm{deg}(Z). 
$$

On the other hand, let us calculate the $L$-function of  $\mathcal M$ modulo $\mathbf m$ using 
Lefschetz trace formula 
$$
 L(C/k, \mathcal M; t) = \prod_i \mathrm{det}(1 - Ft; H_{\mathrm{rig}}^i(C/K, \mathcal M))^{(-1)^{i+1}}
$$
for rigid cohomology \cite[Th\'eor\`eme 6.3]{ELS}. Let $\mathcal N$ be an irreducible subquotient of $\mathcal M$ 
as convergent $F$-isocrystals on $C/K$ such that $\mathcal N$ 
includes the generic slope $0$ part. Since any slope of Frobenius $F_x$ on $i_x^\ast\mathcal N$ 
at any point $x$ of $Z \ne \emptyset$ is positive, $\mathcal N$ is not constant and is of rank $\geq 2$. 
The Poincar\'e duality \cite[Th\'eor\`eme 2.4]{Be2} implies $H^i_{\mathrm{rig}}(C/K,\mathcal N) = 0$ for $i = 0, 2$ 
and then $\mathrm{dim}_K\, H^1_{\mathrm{rig}}(C/K,\mathcal N) \leq 2(g-1)\mathrm{rank}\, \mathcal N$ 
by Euler-Poincar\'e formula of rigid cohomology of curves \cite[Corollaire 5.0-12]{CM}. 
If $\mathcal N' \ne \mathcal N$ is a subquotient of $\mathcal M$, then 
any slope of Frobenius $F_x$ on $i_x^\ast\mathcal N'$ 
at any point $x$ is positive since the generic slope $0$ subquotient of rank $1$ is included in $\mathcal N$. 
It implies that any slope of Frobenius $F$ on the rigid cohomology $H^i_{\mathrm{rig}}(C/K,\mathcal N')$ 
is positive for $i = 0, 1, 2$ if it does not vanish. 
Hence, we have a congruence 
$$
 L(C/k, \mathcal M; t)\,  \equiv\, \mathrm{det}(1 - Ft; H_{\mathrm{rig}}^1(C/K, \mathcal N))\, \, 
 (\mbox{mod}\, \, \mathbf m). 
$$
Since $\mathrm{dim}_K\, H^1_{\mathrm{rig}}(C/K,\mathcal N) \leq 2(g-1)\mathrm{rank}\, \mathcal N \leq 2(g-1)\mathrm{rank}\, \mathcal M$, 
we have a desired inequality. 
}

\begin{corollary}\label{estmax} With the hypothesis of Proposition \ref{est}, $\mathrm{deg}(Z)$ 
is at most $1 + 2(g-1)\mathrm{rank}\mathcal M$. 
\end{corollary} 

\begin{theorem}\label{estjump} Let $C$ be a projective smooth and geometrically connected curve 
over $\mathrm{Spec}\, k$. 
For a convergent $F$-isocrystal $\mathcal M$ on $C/K$, we put a closed subscheme 
$$
      D_{\mathcal M} = \{ x \in C\, |\, \mathrm{NP}(\mathcal M, x) \ne \mathrm{NP}(\mathcal M, \eta)\}, 
$$
regarded as a reduced divisor. Here $\eta$ is the generic point of $C$. Then, for any positive integer $r$, 
there exists a constant $B$ depending on $q$, the genus of $C$ and $r$ such that the inequality 
$$
     \mathrm{deg}(D_{\mathcal M}) \leq B 
$$
holds for any convergent $F$-isocrystal $\mathcal M$ on $C/K$ of rank $r$. 
\end{theorem}

In order to apply Proposition \ref{est} we prepare the lemma below. 

\begin{lemma}\label{genericrank1} Let $Y$ be a geometrically connected scheme separated of finite type 
over $\mathrm{Spec}\, k$ , 
and let $\mathcal N$ be a convergent $F$-isocrystal on $Y/K$. 
Let $\gamma_1 < \cdots < \gamma_2 < \cdots <\gamma_m$ be the generic slopes of $\mathcal N$ 
with multiplicities $r_1, r_2, \cdots, r_m$. 
\begin{enumerate}
\item For $l = 1, \cdots, m$, the exterior power $\wedge^{r_1 + \cdots +r_l} \mathcal N$ of $\mathcal N$ 
has an initial slope of $r_1\gamma_1 + \cdots + r_l\gamma_l$ with multiplicity $1$. 
\item Suppose the Newton polygon function $\mathrm{NP}(\mathcal N, -)$ is not constant on $Y$. 
Then, for some integer $1 \leq l < m$, there exists a point $y$ of $Y$ such that 
the initial slope of $\wedge^{r_1 + \cdots +r_l} \mathcal N$ at $y$ is greater than $r_1\gamma_1 + \cdots + r_l\gamma_l$. 
\end{enumerate}
\end{lemma}

\noindent
{\sc Proof of Theorem \ref{estjump}.} Suppose $\mathcal M$ is an $F$-isocrystal on $C/K$ with 
nonconstant Newton polygons, 
and $\gamma_1 < \cdots < \gamma_2 < \cdots <\gamma_m$ the generic slopes of $\mathcal M$ 
with multiplicities $r_1, r_2, \cdots, r_m$. Then 
$$
   D_{\mathcal M} = \underset{l}{\bigcup} \left\{ x \in C\, \left|\, 
   \begin{array}{l} \mbox{The initial slope of $\wedge^{r_1 + \cdots +r_l} \mathcal M$ 
   at $x$} \\ \mbox{is greater than that at $\eta$.}\end{array} \right. \right\}. 
$$
Applying Lemma \ref{genericrank1} for each $\wedge^{r_1 + \cdots +r_l} \mathcal M$, we 
may suppose the initial generic slope of $\mathcal M$ is $0$ of multiplicity $1$ by a suitable twist of Frobenius. 
Replacing $\mathcal M$ by 
the tensor product $\mathcal M^{\otimes q-1}$, we may assume that 
$\mathcal M$ satisfies the conditions (i) and (ii) of Proposition \ref{est}. 
Then $\mathrm{deg}(Z)$ is bounded depending only on $q$, $g$ and $\mathrm{rank}\, \mathcal M$ 
by Corollary \ref{estmax}. 
Hence, the bound exists since the number of generic slopes is 
less than or equal to $\mathrm{rank}\, \mathcal M$
\hspace*{\fill} $\Box$

\begin{remark}\label{preest}
\begin{enumerate}
\item
More precisely, the upper bounds depend on the cardinal of the base field of definition of $\mathcal M$ 
on $C/K$. If there are a subfield $k_0$ of $k$, a projective smooth curve $C_0$ over $\mathrm{Spec}\, k_0$ and 
a convergent $F$-isocrystal $\mathcal M_0$ on $C_0/K_0$ ($K_0$ is a complete valuation subfield of $K$ with 
the residue field $k_0$) such that $\mathcal M$ is the inverse image of $\mathcal M_0$, 
then $\mathrm{deg}(Z)$ in Corollary \ref{estmax} and $\mathrm{deg}(D_{\mathcal M})$ in Theorem \ref{estjump}
are bounded by the constant depending not on the cardinal of $k$ but on the cardinal of $k_0$. 
\item One can take a upper bound $B = r + 2^{1 + (q-1)r}(g-1)$ 
in Theorem \ref{estjump}. 
This upper bound is not sharp. 
\end{enumerate}
\end{remark}

\subsection{The case of elliptic curves} 
To clarify the idea of proof, we first prove the constancy theorem in the case of elliptic curves. 

\begin{theorem}\label{ECconst} Let $k$ be a finite field, and $X$ an elliptic curve  over $\mathrm{Spec}\, k$. 
Then any convergent $F$-isocrystal $\mathcal M$ on $X/K$ has constant Newton polygons. 
\end{theorem}

\prf{Suppose $\mathcal M$ is an $F$-isocrystal on $X/K$ with 
nonconstant Newton polygons. 
Let $D_{\mathcal M} \ne 0$ be the reduced divisor of $X$ consisting of the points $x$ such that 
$$
     \mathrm{NP}(X, x) \ne \mathrm{NP}(X, \eta)
$$
where $\eta$ is the generic point of $X$. 
Let $[n] : X \rightarrow X$ be the morphism of multiplication $n$ 
for a positive integer $n$. 
If $n$ is prime to $p$, then $[n]$ is a finite etale morphism and we have 
$$
     \mathrm{deg}(D_{[n]^\ast\mathcal M}) = n^2\mathrm{deg}(D_{\mathcal M}). 
$$
It is a contradiction to the boundedness of $\mathrm{deg}(D_{[n]^\ast\mathcal M})$ 
by Theorem \ref{estjump}. 
Therefore, any $F$-isocrystal on $X/K$ has constant Newton polygons. 
}

\subsection{The case of Abelian varieties}

In this section we prove our main theorem for general Abelian varieties.

\begin{theorem}\label{abconst} Let $k$ be a finite field, and $X$ an Abelian variety over $\mathrm{Spec}\, k$ of dimension $g$. 
Any convergent $F$-isocrystal $\mathcal M$ on $X/K$ has constant Newton polygons, i.e., 
$$
  F\mbox{-}\mathrm{Isoc}(X/K)^{\mathrm{CNP}} = F\mbox{-}\mathrm{Isoc}(X/K). 
$$
\end{theorem}

\noindent
{\sc Proof.} Suppose $\mathcal M$ is an $F$-isocrystal on $X/K$ with 
nonconstant Newton polygons. 
Let $D_{\mathcal M} \ne 0$ be the reduced divisor of $X$ consisting of the points $x$ such that 
$$
     \mathrm{NP}(X, x) \ne \mathrm{NP}(X, \eta)
$$
as in the proof of Theorem \ref{ECconst}. Then 
$D_{\mathcal M}$ is purely of codimension $1$ by de Jong-Oort's purity theorem 
(Theorem \ref{Jort}). 
After replacing $k$ by a finite extension, 
we can find a projective smooth geometrically connected curve 
$C$ over $\mathrm{Spec}\, k$ such that 
\begin{enumerate}
\item[(a)] For any closed integral subscheme $Z$ over $\mathrm{Spec}\, k$ of codimension $1$ in $X$, 
the set theoretical intersection $C \cap Z$ is nonempty; 
\item[(b)] $O \in C$ and $O \not\in D_{\mathcal M}$, 
\end{enumerate}
where $O$ is the origin of the Abelian variety $X$. 
Indeed, if we fix an embedding $X$ into a projective space, then 
by Bertini's theorem and Moishezon-Nakai's criterion of ampleness \cite[Appendix A, Theorem 5.1]{Ha},   
one can obtain a projective smooth geometrically connected curve $C$ in $X$ 
as an intersection of different $g-1$ hyperplane sections  
such that the condition (a) holds 
after replacing $k$ by a finite extension. 
Note that the genus of $C$ is greater than or equal to $g$ by weak Lefschetz theorem. 
Since $\mathcal D_{\mathcal M}$ is of codimension $1$, one can choose $C \not\subset D_{\mathcal M}$. 
Taking a translation by a $k$-rational point after replacing $k$ by a finite extension if necessary, 
we may also assume the condition (b). 

\begin{lemma}\label{cont} 
\begin{enumerate}
\item For any positive integer $n$, we have $C \cap D_{[n]^\ast\mathcal M} \ne \emptyset$ 
and $C \not\subset D_{[n]^\ast\mathcal M}$. 
\item If we put $\Omega = \{ x \in C(\overline{k})\, |\, x \in C \cap D_{[n]^\ast\mathcal M}\, \, \mbox{for some}\, \, n > 0\}$, 
then $\Omega$ is infinite. 
\item If $D_{[n]^\ast\mathcal M} \cap C$ denotes the set theoretical intersection 
and $(D_{[n]^\ast\mathcal M} \cap C)_{\mathrm{red}}$ denotes a reduced divisor of $C$, then 
$$
   \underset{n}{\mathrm{sup}}\, \mathrm{deg}((D_{[n]^\ast\mathcal M} \cap C)_{\mathrm{red}}) = \infty. 
$$
\end{enumerate}
\end{lemma}

\prf{(1) It follows from the condition (a) and $O \not\in D_{[n]^\ast\mathcal M} \cap C$ for any $n$. 

(2) Suppose that $\Omega$ is finite, namely, $\Omega = \{ y_1, y_2, \cdots, y_s\}$. 
Let $n$ be a positive integer such that $n$ is a multiple of orders of all $y_1, y_2, \cdots, y_s$. 
Since $[n](y_i) = O$, $y_i \not\in D_{[n]^\ast\mathcal M} \cap C$ for any $i$. 
This contradicts to $C \cap D_{[n]^\ast\mathcal M} \ne \emptyset$. 

(3) For a finite extension $k_r$ of $k$ of degree $r$, the cardinal of $X(k_r)$ is finite. Hence 
there is a properly infinitely increasing sequence $r_1 < r_2 < \cdots$ such that 
$x_{r_i} \in \Omega \cap X(k_{r_i})$ and $x_{r_i} \not\in X(k')$ 
for any proper subfield $k'$ of $k_{r_i}$ by (2). Since 
$C$ and $D_{[n]^\ast\mathcal M}$ is defined over $\mathrm{Spec}\, k$, 
the cardinal of the set theoretical intersection $D_{[n]^\ast\mathcal M} \cap C$ is greater than or equal to $r_i$ 
if $x_{r_i} \in D_{[n]^\ast\mathcal M} \cap C$ because of Galois conjugation. 
}

\vspace*{3mm}

Now we return to the proof of Theorem \ref{abconst}. 
Since $\mathrm{rank}\, [n]^\ast\mathcal M = \mathrm{rank}\, \mathcal M$ for any positive integer $n$, the 
degree of $(D_{[n]^\ast\mathcal M} \cap C)_{\mathrm{red}}$ in $C$ is bounded 
by a constant independent of $n$ by Theorem \ref{estjump} and Remark \ref{preest}. 
This contradicts to Lemma \ref{cont} (3). 
Therefore, any $F$-isocrystal on $X/K$ has constant Newton polygons. 
\hspace*{\fill} $\Box$

\section{Isotriviality of a family of curves} 

In this section we study the isotriviality of projective smooth families of connected  
curves on Abelian varieties of characteristic $p > 0$. 
Let $k$ be an algebraically closed field of characteristic $p$. 

\subsection{Isotriviality of families of ordinary abelian varieties}

We recall some definitions.

\begin{definition}\label{ordi} 
Let $\kappa$ be a perfect field of characteristic $p$ and $\overline{\kappa}$ an algebraic closure of $\kappa$. 
\begin{enumerate}
\item An Abelian variety $S$ of dimension $g$ over $\mathrm{Spec}\, \kappa$ 
is said to be ordinary if $\mathrm{dim}_{\mathbb F_p}S[p](\overline{\kappa}) = g$. 
Here $S[p]$ is the subgroup scheme defined by the kernel of multiplication with $p$. 
\item A projective smooth geometrically connected curve over $\mathrm{Spec}\, \kappa$ 
is said to be ordinary if the Jacobian variety $J(S)$ of $S$ is so.
\end{enumerate}
\end{definition}

\begin{theorem}\label{isNP}  Let $S$ be 
a projective smooth connected scheme over $\mathrm{Spec}\, k$ such that 
any geometric convergent $F$-isocrystal on $S/K$ has constant Newton polygons. 
Suppose $f : X \rightarrow S$ is a polarized Abelian scheme relatively 
of dimension $g$ and of degree $d$.  
If there is a point $s \in S$ such that the fiber $X_s$ of $X$ at $s$ is an ordinary Abelian variety, 
then $X$ is isotrivial over $S$. 
\end{theorem}

\prf{Let us consider the relative first rigid cohomology $R^1f_{\mathrm{rig} \ast}\mathcal O_{]X[}$ 
of the morphism $f : X \rightarrow S$ which is a convergent 
$F$-isocrystal on $S/K$ (see Theorem \ref{coh}). Since there is a point $s \in S$ such that $X_s$ is ordinary, 
$X$ is a family of ordinary Abelian varieties by our hypothesis on Newton polygons. 
Because the ordinary locus $\mathcal A_{g, d, k}^{\mathrm{ord}}$ of the coarse moduli scheme 
$\mathcal A_{g, d, k}$ of 
polarized Abelian varieties of dimension $g$ and of degree $d$ 
over $\mathrm{Spec}\, k$ is quasi-affine by \cite[XI, Th\'eor\`eme 5.2]{MB}, 
the canonical morphism $S \rightarrow \mathcal A_{g, d, k}^{\mathrm{ord}}$ is constant.}

\begin{corollary}\label{isJNP}  Let $S$ be 
a projective smooth connected scheme of finite type over $\mathrm{Spec}\, k$ such that 
any geometric convergent $F$-isocrystal on $S/K$ has constant Newton polygons. 
Suppose $X \rightarrow S$ is a projective smooth family of curves of genus $g$.  
If there is a point $s \in S$ such that the fiber $X_s$ of $X$ at $s$ is ordinary, 
then $X$ is isotrivial over $S$. 
\end{corollary}

\prf{Consider the relative Jacobian variety $J(X/S)$ of $X$ over $S$. Then $J(X/S)$ is a family of polarized Abelian scheme 
over $S$ with a ordinary fiber $J(X/S)_s$. Then the family $J(X/S)$ is isotrivial by Proposition \ref{isNP}. 
Now apply the Torelli theorem (see \cite[Theorem 12.1]{Mi}, for example), and then the family $X$ is isotrivial over $S$. 
}

\begin{corollary}\label{isocNP} Let $S$ be 
an Abelian variety over $\mathrm{Spec}\, k$, and 
$X$ a polarized Abelian scheme over $S$. 
If there exists a point $s \in S$ such that the fiber $X_s$ of $X$ at $s$ is an ordinary Abelian variety, then 
$X$ is isotrivial over $S$. 
\end{corollary}

\prf{We have only to prove any geometric convergent $F$-isocrystal on $S/K$ has 
constant Newton polygons by Theorem \ref{isNP}. Let $Y \rightarrow S$ be a proper smooth morphism. 
Then there exist a smooth integral scheme $\mathcal T$ separated of finite type over a spectrum 
$\mathrm{Spec}\, \kappa$ 
of a finite field $\kappa$ such that 
$\kappa(\mathcal T) \subset k$, an Abelian scheme $\mathcal S$ over $\mathcal T$ and 
a proper smooth morphism $g : \mathcal Y \rightarrow \mathcal S$ 
such that the base change by the morphism $\mathrm{Spec}\, k \rightarrow \mathcal T$ is 
the given morphism $Y \rightarrow S$. Here $\kappa(\mathcal T)$ is the function field of $\mathcal T$. 
Then, for any closed point $t \in \mathcal T$,  
the relative rigid cohomology $\mathcal M = R^ig_{\mathrm{rig}\ast}\mathcal O_{]\mathcal Y[}$ 
has a constant Newton polygon on $\mathcal S_t$ by Theorem \ref{abconst}. 
Hence, any geometric convergent $F$-isocrystal on $S/K$ 
has constant Newton polygons by Proposition \ref{geom}. 
}

\subsection{Isotriviality of family of curves over Abelian varieties}

We state the isotriviality of a family of curves. In the case $g \leq 1$ the isotriviality 
below is trivial without the hypotheses (i) and (ii) since the coarse moduli scheme $\mathfrak M_g$ 
of projective smooth curves of genus $g$ is affine. 
We will prove the following theorem in the rest of this section. 

\begin{theorem}\label{isotcNP} Let $S$ be 
a projective smooth and connected scheme over $\mathrm{Spec}\, k$, 
and $f : X \rightarrow S$ a projective smooth family of connected curves of genus $g \geq 2$. 
Suppose that any geometric convergent $F$-isocrystal on $S/K$ 
has constant Newton polygons. 
Then the family $X$ over $S$ is isotrivial. 
\end{theorem}

\begin{corollary}\label{isotab} Let $S$ be an Abelian variety over $\mathrm{Spec}\, k$.
If $f : X \rightarrow S$ is a projective smooth family of connected curves, 
then the family $X$ over $S$ is isotrivial. 
\end{corollary}

\prf{Since any geometric convergent $F$-isocrystal on $S/K$ has 
constant Newton polygons by the proof of Corollary \ref{isocNP}, 
the assertion follows from Theorem \ref{isotcNP}. 
}

\subsection{Torelli theorem} 
In this subsection we recall Tamagawa's work in \cite{tam}. 
Let $X$ be a projective smooth curve of genus $g \geq 2$ 
over $\mathrm{Spec}\, k$. 

\begin{definition}\label{gon} The gonality of $X$ is a minimum degree of 
nonconstant morphisms $X \rightarrow \mathbb P^1_k$. 
\end{definition}

\begin{theorem}\label{exon} \mbox{\rm (\cite[Proposition 2.14]{tam})} Suppose that 
the $p$-rank of $X$ is neither $0$ nor $g$. For any sufficient large integer $n$ which is prime to $p$, 
there exists a finite etale morphism $Y \rightarrow X$ cyclic of degree $n$ 
such that the gonality of $Y$ is greater than or equal to $5$. 
\end{theorem}

\begin{remark} In \cite[Proposition 2.14]{tam} the existence of $Y$ with gonality $\geq 5$ 
is proved with the assumption that $X$ is non-almost elliptic, i.e., the Jacobian variety 
$J(X)$ is not isogenous to $E^g$ for an elliptic curve $E$ over $\mathrm{Spec}\, k$. 
\end{remark}

\begin{definition}\label{newordtor} Let  $f : Y \rightarrow X$ be a finite etale morphism 
of projective smooth and connected curves over $\mathrm{Spec}\, k$, 
$J(X) = \mathrm{Pic}^0(X)$ and 
$J(Y) = \mathrm{Pic}^0(Y)$ the Jacobian varieties of $X$ and $Y$, respectively, and define 
the new part of Jacobian variety of $Y$ relatively to $X$ by the Abelian variety 
$$
       J(Y, X) = J(Y)/f^\ast J(X)  
$$
with a polarization. 
The finite etale covering $Y$ over $X$ is said to be new-ordinary if $J(Y, X)$ is an ordinary Abelian variety. 
\end{definition}

\begin{theorem}\label{neword} \mbox{\rm (\cite[Theorem 4.3.1]{ray}, \cite[Corollary 5.3]{tam})} 
Let $X$ be a projective smooth curve of genus $g \geq 2$. 
For any sufficiently large prime number $l$ which is prime to $p$, there exists a nontrivial 
$\mu_l$-torsor $Y$ of $X$ 
which is new-ordinary.  
\end{theorem} 

Let $\mathcal C_k$ be a category of Artinian local ring with residue field $k$. 
For a proper smooth connected scheme $T_0$, we define a deformation functor 
$$
      M_{T_0} : \mathcal C_k \rightarrow \mbox{(Sets)}, 
$$
that is, for $R \in \mathcal C_k$, the set $M_{T_0}(R)$ is a set of isomorphism classes of pairs $(T, \varphi)$ 
such that $T$ is a proper smooth scheme over $\mathrm{Spec}\, R$ and $\varphi$ is an isomorphism 
$T \times_{\mathrm{Spec}\, R}\mathrm{Spec}\, k \cong T_0$. If $T_0$ is a projective smooth connected curve 
or an Abelian variety, then $M_{T_0}$ is pro-representable by the formal spectrum of a ring of formal power series over $k$. 

Now let $l$ be a prime number which is prime to $p$, $X_0$ a projective smooth connected curve of genus $\geq 2$ 
over $\mathrm{Spec}\, k$, $\mu_l$ the locally constant etale sheaf consisting of $l$-th roots of unity, 
$Y_0$ a $\mu_l$-torsor of $X_0$ which associates to a nontrivial element $\mathcal L_0 \in J(X_0)[l](k) 
\cong H^1_{\mathrm{et}}(X_0, \mu_l)$, and $J(Y_0, X_0)$ the new part of Jacobian variety 
of $Y_0$ with respect to $X_0$. Then there is a canonical map 
$$
       T_{\mathcal L_0}(R) :  M_{X_0}(R) \rightarrow M_{J(Y_0, X_0)}(R)
$$
as follows. For $X \in M_{X_0}(R)$, there exists a unique $\mu_l$-torsor $f : Y \rightarrow X$ over $\mathrm{Spec}\, R$
which is a lift of the $\mu_l$-torsor $f_0 : Y_0 \rightarrow X_0$ over $\mathrm{Spec}\, k$ 
by \cite[I, Corollaire 8.4]{SGA1}. Then $T_{\mathcal L_0}(R)(X)$ 
is the new part $J(Y, X) = J(Y)/f^\ast J(X)$ of the relative Jacobian variety 
$J(Y)$ of $Y$ with respect to the relative Jacobian variety $J(X)$ of $X$ with a natural isomorphism 
$J(Y, X) \times_{\mathrm{Spec}\, R}\mathrm{Spec}\, k \cong J(Y_0, X_0)$. 
Hence one has a Torelli morphism 
$$
         T_{\mathcal L_0} : M_{X_0} \rightarrow M_{J(Y_0, X_0)}. 
$$
By the numerical estimates in \cite[Corollaries 4.7, 5.3]{tam} one has 

\begin{theorem}\label{tore} \mbox{\rm (\cite[Theorem 3.3]{saidi})} With the notation as above, 
suppose that the gonality of $X_0$ is greater than or equal to $5$. 
Then, for any sufficiently large prime number $l$ which is prime to $p$, 
there exists a nontrivial element $\mathcal L_0 \in J(X_0)[l](k)$ such that 
the corresponding $\mu_l$-torsor $Y_0$ of $X_0$ is new-ordinary and 
the Torelli morphism $T_{\mathcal L_0}$ is an immersion. 
\end{theorem} 

\subsection{Proof of Theorem \ref{isotcNP}.} 
First we prepare several assertions. 
Lemmas \ref{sai2} and \ref{etind} are higher dimensional generalizations 
of Lemma/Proposition 4.7 and Lemma 2 in \cite{saidi}, respectively. 

\begin{lemma}\label{sai2} Let $S$ be a smooth connected 
scheme separated of finite type over $\mathrm{Spec}\, k$, 
$f : X \rightarrow S$ be a projective smooth morphism over $\mathrm{Spec}\, k$ such that 
each fiber of $f$ is a connected curve, $t \in X(k(S)^{\mathrm{sep}})$ a point of degree $d_t$ over $k(S)$ 
where $k(S)^{\mathrm{sep}}$ is a separable closure of the function field $k(S)$ of $S$, 
and $n$ a positive integer which is prime to $pd_t$.   
Let us fix a closed point $s$ of $S$ and a $\mu_n$-torsor $Y_s$ of the fiber $X_s$ of $f$ at $s$. 
Then there exist a finite etale morphism $S' \rightarrow S$ with a base change morphism 
$f' : X' = X \times_S S' \rightarrow S'$ and a $\mu_n$-torsor 
$Y' \rightarrow X'$ such that, for $s' \in S'$ which goes to $s$ in $S$, the $\mu_n$-torsor $Y'_{s'}$ 
of $X'_{s'}$ is naturally isomorphic to $Y_s$ over $X_s$. 
\end{lemma}

\prf{Since we may identify $\mu_n$ with the constant etale sheaf $\mathbb Z/n\mathbb Z$, 
$R^1f_{\mathrm{et}\ast}\mu_n$ is locally constant. Moreover, we can take 
a finite etale cover $S'$ of $S$ such that $R^1f_{\mathrm{et}\ast}\mu_n$ is constant. 
Take a point $s' \in S'$ with a closed immersion $i_{s'} : s' \rightarrow S'$ 
such that $s'$ is above $s$. Then we have an 
isomorphism $i_{s'}^\ast R^1f_{\mathrm{et}\ast}'\mu_n \cong H_{\mathrm{et}}^1(X_{s'}, \mu_n)$ 
by proper base change theorem of etale cohomology. 
Since $R^1(f')_{\mathrm{et}\ast}\mu_n$ is constant, the natural homomorphism 
$$
      H_{\mathrm{et}}^0(S', R^1f_{\mathrm{et}\ast}'\mu_n) \rightarrow 
      H_{\mathrm{et}}^1(X_{s'}', \mu_n) \cong H_{\mathrm{et}}^1(X_s, \mu_n)
$$
is bijective. By compositing with the exact sequence 
$$
    0 \rightarrow H^1_{\mathrm{et}}(S', \mu_n) \rightarrow H^1_{\mathrm{et}}(X', \mu_n) 
    \rightarrow H_{\mathrm{et}}^0(S', R^1f_{\mathrm{et}\ast}'\mu_n) 
    \rightarrow H^2_{\mathrm{et}}(S', \mu_n) \rightarrow H^2_{\mathrm{et}}(X', \mu_n) 
$$
arising from Leray spectral sequence, 
we have only to prove that 
the natural homomorphism $H^2_{\mathrm{et}}(S', \mu_n) \rightarrow H^2_{\mathrm{et}}(X', \mu_n)$ is injective. 
Then the homomorphism 
$$
      H^1_{\mathrm{et}}(X', \mu_n) \rightarrow H_{\mathrm{et}}^1(X_s, \mu_n)
$$
is surjective, and hence there exists a $\mu_n$-torsor $Y' \rightarrow X'$ whose fiber at $s'$ is isomorphic to 
the given $\mu_n$-torsor $Y_s \rightarrow X_s$. 

Let $T$ be the normalization of $S$ in the function field $k(S)(t)$ of $t$. 
Then the restriction $g = f|_T : T \rightarrow S$ is generically etale and finite of degree $d_t$. 
If we put $T' = T\times_S S'$ with a finite morphism $g' : T' \rightarrow S'$, then we have a natural commutative diagram 
$$
     \begin{array}{ccccc}
          H^2_{\mathrm{et}}(S', \mu_n) &\overset{(f')^\ast}{\rightarrow} &H^2_{\mathrm{et}}(X', \mu_n) \\
          &(g')^\ast \searrow \hspace*{5mm}&\downarrow \\
          & &H^2_{\mathrm{et}}(T', \mu_n). 
     \end{array}
$$
The injectivity of $(f')^\ast$ follows from that of $(g')^\ast$ which is proved in Proposition \ref{inj} below. 
}

\begin{proposition}\label{inj} 
Let $S$ be a smooth connected scheme separated of finite type over the 
spectrum $\mathrm{Spec}\, k$ of an algebraic closed field $k$ of characteristic $p \geq 0$, 
and $g : T \rightarrow S$ a generically etale and finite morphism of degree $d$ 
such that $T$ is normal. 
Suppose $n$ is a positive integer which is prime to $dp$ (resp. $d$) if $p > 0$ (resp. $p = 0$). 
Then the homomorphism 
$$
       g^\ast : H_{\mathrm{et}}^2(S, \mu_n) \rightarrow H_{\mathrm{et}}^2(T, \mu_n)
$$
is injective. 
\end{proposition}

\prf{
Let $U$ be an open dense subscheme of $S$ such that the inverse image 
$g^{-1}(U)$ is regular and the complement $E$ of $U$ in $S$ 
is of codimension $\geq 2$. Such a $U$ exists since $T$ is normal. 
Then we have a natural commutative diagram 
$$
     \begin{array}{ccccc}
          H^2_{\mathrm{et}}(S, \mu_n) &\overset{g^\ast}{\rightarrow} &H^2_{\mathrm{et}}(T, \mu_n) \\
          \downarrow & &\downarrow \\
          H^2_{\mathrm{et}}(U, \mu_n) &\underset{g^\ast}{\rightarrow} &H^2_{\mathrm{et}}(g^{-1}(U), \mu_n). 
     \end{array}
$$
Since the right vertical homomorphism is injective by the cohomological purity theorem \cite[VI, Theorem 5.1]{mil}, 
we may assume that $T$ is smooth over $\mathrm{Spec}\, k$. 

Let us take a largest open dense subscheme $V$ of $S$ such that 
the inverse image $W = g^{-1}(V)$ is etale over $S$, and 
$E$ (resp. $F$) is the complement of $V$ (resp. $W$) in $S$ (resp. $T$). 
Let us now consider the commutative diagram 
$$
\begin{array}{ccccc}
H_{\mathrm{et}, E}^2(S, \mu_n) &\rightarrow &H_{\mathrm{et}}^2(S, \mu_n)
&\rightarrow &H_{\mathrm{et}}^2(V, \mu_n) \\
g^\ast\downarrow\hspace*{2mm} & &g^\ast\downarrow\hspace*{2mm}  & &\hspace*{2mm} \downarrow g^\ast \\
H_{\mathrm{et}, F}^2(T, \mu_n) &\rightarrow &H_{\mathrm{et}}^2(T, \mu_n)
&\rightarrow &H_{\mathrm{et}}^2(W, \mu_n) \\
g_\ast\downarrow\hspace*{2mm} & &g_\ast\downarrow\hspace*{2mm}  & &\hspace*{2mm} \downarrow g_\ast \\
H_{\mathrm{et}, E}^2(S, \mu_n) &\rightarrow &H_{\mathrm{et}}^2(S, \mu_n)
&\rightarrow &H_{\mathrm{et}}^2(V, \mu_n), 
\end{array}
$$
where the left items are etale cohomology with supports, each sequence of horizontal homomorphisms 
is a localization sequence, and the homomorphisms $g_\ast : H_{\mathrm{et}, c}^{2\mathrm{dim}T}(T, \mu_n) 
\rightarrow H_{\mathrm{et}, c}^{2\mathrm{dim}S}(S, \mu_n)$ and so on are induced by 
the pullbacks 
$g^\ast : H_{\mathrm{et}, c}^{2\mathrm{dim}S}(S, \mu_n) \rightarrow H_{\mathrm{et}, c}^{2\mathrm{dim}T}(T, \mu_n)$ and so on 
of etale cohomology with compact supports under Poincar\'e duality 
and trace morphisms \cite[VI, Remark 11.6]{mil}. 
Note that Poincar\'e duality is applicable since $g : T \rightarrow S$ is a finite morphism of smooth schemes. 
Then the composite $g_\ast g^\ast$ in the left (resp. right) vertical homomorphisms 
is the map of multiplication with $d$ by Lemma \ref{sure} below (resp. 
the finite etaleness of $W$ over $V$). 
Hence, the middle $g_\ast g^\ast$ is surjective. Indeed, for any $a \in H_{\mathrm{et}}^2(S, \mu_n)$, there exists 
an element $b \in H_{\mathrm{et}}^2(S, \mu_n)$ coming from $H_{\mathrm{et}, E}^2(S, \mu_n)$ 
such that $g_\ast g^\ast(a + b) = da$. Therefore the finiteness of $H^2_{\mathrm{et}}(S, \mu_n)$ 
implies that 
the homomorphism $g^\ast : H^2_{\mathrm{et}}(S, \mu_n) \rightarrow H^2_{\mathrm{et}}(T, \mu_n)$ 
is injective.
}

\begin{lemma}\label{sure} With the notation in Proposition \ref{inj}, the followings hold.  
\begin{enumerate}
\item Any generic point of complement of $V$ in $S$ is pure of codimension $1$. 
The same holds for the complement of $W$ in $T$. 
\item Let $E_1, \cdots, E_r$ be reduced irreducible components of the complement $E$ of $V$ in $S$, 
and $F_{i,1}, \cdots, F_{i,s_i}$ the reduced irreducile components of the inverse image $F_i$ of $E_i$ in  $T$ 
with multiplicity $e_{i, j}$ of $F_{i, j}$ in $T \times_S E_i$ and the degree of $F_{i, j}$ over $E_i$ (here we consider 
the ramification index and the residual degree for the extension $\mathcal O_{T, F_{i, j}}/\mathcal O_{S, E_i}$ 
of discrete valuation rings, respectively). Then $\sum_j e_{i, j}f_{i, j} = d$ for any $i$.
\item The homomorphism 
$g^\ast : H_{\mathrm{et}, E_i}^2(S, \mu_n) \rightarrow H_{\mathrm{et}, F_i}^2(T, \mu_n)$ 
is given by the homomorphism 
$$
       g^\ast : H_{\mathrm{et}}^0(E_i, \mu_n) \rightarrow \oplus_j\, H_{\mathrm{et}}^0(F_{i, j}, \mu_n) 
$$
under the isomorphism (and the same for $F_{i, j}$) induced by the bottom horizontal 
Gysin isomorphism \cite[VI, Theorem 5.1]{mil}:
$$
      \begin{array}{ccccc}
          \mu_n(k) &\cong &H_{\mathrm{et}}^0(E_i, \mu_n)(-1) &\rightarrow &H_{\mathrm{et}, E_i}^2(S, \mu_n) \\
          &&\cong \downarrow \hspace*{3mm} & &\hspace*{3mm} \downarrow \cong \\
          &&H_{\mathrm{et}}^0(E_i^{\mathrm{sm}}, \mu_n)(-1)&\overset{\cong}{\rightarrow} 
          &H_{\mathrm{et}, E_i^{\mathrm{sm}}}^2(S \setminus E_i^{\mathrm{sing}}, \mu_n), 
      \end{array}
$$
where $E_i^{\mathrm{sing}}$ is a singular locus of $E_i$ (note that $E_i^{\mathrm{sing}}$ is of codimension $\geq 2$ in $S$), 
$E_i^{\mathrm{sm}} = E_i \setminus E_i^{\mathrm{sing}}$, 
and $(-1)$ means the $(-1)$-th Tate twist. 
\item The homomorphism 
$g_\ast : H_{\mathrm{et}, F_i}^2(T, \mu_n) \rightarrow H_{\mathrm{et}, E_i}^2(S, \mu_n)$ 
is given by the homomorphism 
$$
       \oplus_jH_{\mathrm{et}}^0(F_{i, j}, \mu_n) \rightarrow H_{\mathrm{et}}^0(E_i, \mu_n) 
       \hspace*{6mm} (a_i) \mapsto \sum_i a_i 
$$
under the isomorphisms as in (3). 
\item The composite $g_\ast g^\ast : H_{\mathrm{et}, E}^2(S, \mu_n) \rightarrow 
H_{\mathrm{et}, E}^2(S, \mu_n)$ is the map of multiplication with $d$. 
\end{enumerate}
\end{lemma} 

\prf{(1) Since $U$ is the largest and $T$ is finite over $S$, 
the assertion follows from Zariski-Nagata purity theorem \cite[X, Th\'eor\`eme 3.4]{SGA2}. 

(4) The homomorphism is induced by the natural homomorphism 
$$
        g^\ast : H_{\mathrm{et}, c}^{2\mathrm{dim} E_i}(E_i, \mu_n) \rightarrow 
        H_{\mathrm{et}, c}^{2\mathrm{dim} F_i}(F_i, \mu_n) \cong \oplus_j H_{\mathrm{et}, c}^{2\mathrm{dim} F_{i, j}}(F_{i, j}, \mu_n). 
$$

(5) Since $H_{\mathrm{et}, E}^2(S, \mu_n) \cong \oplus_i, H_{\mathrm{et}, E_i}^2(S, \mu_n)$, 
the assertion follows from (2), (3), and (4). 
}

\begin{lemma}\label{etind} Let $S$ be a projective smooth scheme over $\mathrm{Spec}\, k$, 
and $X$ a proper smooth family of curves of genus $g \geq 2$. 
Let $S' \rightarrow S$ be a finite etale morphism and $Y' \rightarrow X' = X \times_SS'$ 
a finite etale morphism. If the family $Y'$ over $S'$ is isotrivial, then so is the family $X$ over $S$. 
\end{lemma} 

\prf{When $S$ is a curve, the lemma is just Lemma 2 in \cite[p.235]{saidi}. 
Let $T$ be a projective smooth curve in $S$, and $T', X_T, X_T', Y'_T$ base changes of $S', X, X', Y'$ 
by $T \rightarrow S$, respectively. Then the family $X_T$ over $T$ is isotrivial. 
By varying projective smooth curves in $S$,  there exists an open dense subscheme $U$ of $S$ such that 
$X_U = X \times_SU$ over $U$ is isotrivial by Bertini's theorem. Hence, the family $X$ over $S$ is isotrivial 
because the canonical morphism $S \rightarrow \mathfrak M_g$ is constant. 
Here $\mathfrak M_g$ is the coarse moduli scheme of projective smooth curves of genus $g$. 
}

\vspace*{3mm}

Now we prove Theorem \ref{isotcNP}. 
The following arguments are essentially due to the proof of \cite[Theomre 4.6]{saidi}. 
Let us fix a closed point $s$ in $S$. If the fiber $X_s$ of $X$ at $s$ is ordinary, then 
the assertion follows from Corollary \ref{isJNP}. 
Hence we may suppose the fiber $X_s$ is not ordinary. 
Then there exist 
\begin{enumerate}
\item[$1^\circ$] (only in the case where the $p$-rank of $X_s$ is $0$; 
if not, then $S_0 = S$ and $s_0 = s$) a point $t_0$ of $X(k(S)^{\mathrm{sep}})$ 
of degree $d_0$, a finite etale morphism 
$S_0 \rightarrow S$ and a finite etale morphism 
$Y_0 \rightarrow X \times_SS_0$ such that, for a closed point $s_0$ of $S_0$ which goes to $s$ in $S$, 
$Y_0 \rightarrow X \times_SS_0$ is a $\mu_{l_0}$-torsor 
for a prime number $l_0$ which is prime to $pd_0$ and that 
the $p$-rank of the fiber $Y_{0, s_0}$ of $Y_0$ at $s_0$ is neither $0$ nor $g$;
\item[$2^\circ$] a point $t_1$ of $Y_0(k(S_0)^{\mathrm{sep}})$ 
of degree $d_1$, a finite etale morphism $S_1 \rightarrow S_0$ and a finite etale morphism 
$Y_1 \rightarrow X_1 = Y_0 \times_{S_0}S_1$ such that, for a closed point $s_1$ of $S_1$ which goes to $s_0$ in $S_0$, 
$Y_1 \rightarrow X_1$ is a $\mu_{l_1}$-torsor 
for a prime number $l_1$ which is prime to $pd_1$ and that 
the gonality of the fiber $Y_{1, s_1}$ of $Y_1$ at $s_1$ is greater than or equal to $5$;
\item[$3^\circ$] a point $t_2$ of $Y_1(k(S_1)^{\mathrm{sep}})$ 
of degree $d_2$, a finite etale morphism $S_2 \rightarrow S_1$, a prime number $l_2$ which is prime to $pd_2$, 
and a nontrivial $\mu_{l_2}$-torsor $Y_2 \rightarrow X_2 = Y_1 \times_{S_1}S_2$ such that, 
if $s_2$ is a closed point of $S_2$ which goes to $s_1$ in $S_1$ and 
$\mathcal L \in J(X_{2, s_2})[l_2](k)$ corresponds to 
the $\mu_{l_2}$-torsor $Y_{2, s_2} \rightarrow X_{2,s_2}$, then 
(i) $Y_{2, s_2} \rightarrow X_{2,s_2}$ is new-ordinary and 
(ii) the Torelli morphism $T_{\mathcal L} : M_{X_{2, s_2}} \rightarrow M_{J(Y_{2, s_2}, X_{2,s_2})}$ is an immersion.
\end{enumerate}
by Theorem \ref{neword} and Lemma \ref{sai2} for $1^\circ$, 
by Theorem \ref{exon} and Lemma \ref{sai2} for $2^\circ$ 
and by Theorem \ref{tore} and Lemma \ref{sai2} for $3^\circ$. 
Note that the gonality of $X_{2, s_2} \cong Y_{1, s_1}$ is greater than or equal to $5$ by $2^\circ$. 
In this situation we have only to prove the family $X_2$ over $S_2$ is isotrivial by Lemma \ref{etind}. 
Let $J(X_2)$ and $J(Y_2)$ be relative Jacobian varieties over $S_2$ and define the new part 
$$
          J(Y_2, X_2) = J(Y_2)/f_2^\ast J(X_2)
$$
where $f : Y_2 \rightarrow X_2$ is the canonical morphism. 
Then $J(Y_2, X_2)$ is a polarized Abelian scheme over $S_2$. 
Since $J(Y_2, X_2)_{s_2} = J(Y_{2, s_2}, X_{2,s_2})$ is an ordinary Abelian variety, the family $J(Y_2, X_2)$ over $S_2$ 
is isotrivial by Theorem \ref{isNP}. Hence the family $X_2$ is isotrivial over $S_2$ 
since the Torelli morphism $T_{\mathcal L}$ is an immersion by our construction. 
This completes a proof. 
\hspace*{\fill} $\Box$

\subsection{Existence of convergent $F$-isocrystals with non constant Newton polygons}\label{nonisotsec}

In the study of nonconstant geometric etale fundamental groups on a family, 
M.Sa\"idi proved the following theorem. 

\begin{theorem}\label{nonisot} \mbox{\rm (\cite[Theorem 4.6]{saidi})} Let $C$ 
be a projective smooth and connected curve over $\mathrm{Spec}\, k$, and $f : X \rightarrow C$ 
a proper smooth family of connected curves of genus $\geq 2$. 
If $X$ is not isotrivial over $C$, then there exist a finite etale covering $C'$ of $C$ 
and a proper smooth family $f' : X' \rightarrow C'$ of connected curves 
such that the $p$-ranks of fibers $X_s'\, (s \in C')$ are not constant. 
\end{theorem}

Since there always exists a projective smooth and connected curve $C$ with 
a nonisotrivial family of curves (see \cite[Theorem 3.1]{oort},  \cite{Dia}), 
we have an existence theorem below. 
Moreover, such a curve $C$ and a convergent $F$-isocrystal are possibly defined 
over a finite base field in any characteristic $p$. 

\begin{corollary}\label{extncnp} 
There exist a projective smooth and connected curve $C$ over $\mathrm{Spec}\, k$ 
and a convergent $F$-isocrystal on $C/K$ with nonconstant Newton polygons.
\end{corollary}

\vspace*{5mm}


\appendix

\begin{center}
{\large Appendix}
\end{center}

\section{Proof of Theorem \ref{fil}.}\label{esf}

\subsection{First reduction}\label{filas}

In this appendix we will prove the following theorem. 

\begin{theorem}\label{fila} 
Let $X$ be a smooth scheme separated of finite type over $k$, and $\mathcal M$ 
a convergent $F$-isocrystal on $X/K$. Suppose that 
\begin{enumerate}
\item[($\ast$)] the initial slope of $\mathcal M$ 
at the generic point of $X$ is $0$ of rank $r_0$ 
and the multiplicity of slope $0$ in $F$-isocrystal $i_x^\ast\mathcal M$ is constant on $x \in X$. 
\end{enumerate}
Then there exists a unit-root convergent sub $F$-isocrystal $\mathcal L$ of $\mathcal M$ on $X/K$ of rank $r_0$. 
\end{theorem}

Theorem \ref{fil} follows from the theorem above. Indeed, 
by taking a ramified finite extension $K'$ of $K$ 
with an extension of Frobenius such that the valuation corresponding 
to the initial slope is contained in the valuation group of $K'$, 
we can reduce the assertion to that in the case where the initial slope is $0$ 
by Lemma \ref{sca} below. 
Hence Theorem \ref{fila} above implies Theorem \ref{fil}. 

\begin{lemma}\label{sca} Let $K'$ be a finite extension of $K$ with a residue field $k'$ such that there 
exists a $q$-Frobenius $\sigma'$ on $K'$ satisfying $\sigma'|_K = \sigma$, 
and put $X' = X \times_{\mathrm{Spec}\, k}\mathrm{Spec}\, k$. 
Let $\mathcal M$ be a convergent $F$-isocrystal on $X/K$, and $\mathcal M'$ a convergent $F$-isocrystal 
on $X'/K'$ which is an inverse image of $\mathcal M$. Suppose there exists a convergent sub $F$-isocrystal $\mathcal L'$ 
of $\mathcal M'$ such that all slopes of $i_{\overline{x}}^\ast(\mathcal M'/\mathcal L')$ 
is greater than those of $i_{\overline{x}}^\ast\mathcal L'$ for any geometric point $\overline{x}$ of $X$. 
Then there exists a convergent sub $F$-isocrystal $\mathcal L$ of $\mathcal M$ 
such that the inverse image of $\mathcal L$ on $X'/K'$ is isomorphic to $\mathcal L'$. 
\end{lemma}

\prf{Let $g : X' \rightarrow X$ be a canonical morphism. If we put $\mathcal L$ 
to be the kernel of the natural homomorphism $\mathcal M \rightarrow g_\ast\mathcal M'/g_\ast\mathcal L'$, i.e., 
$$
      \mathcal L = \mathrm{Ker}(\mathcal M \rightarrow g_\ast\mathcal M'/g_\ast\mathcal L'), 
$$
then $\mathcal L$ is the desired convergent sub $F$-isocrystal of $\mathcal M$ by the hypothesis of slopes.}

\vspace*{3mm}

We may also assume that 
\begin{enumerate}
\item[(i)] $\mathbb F_q \subset k$ and $K_\sigma \otimes_{W(\mathbb F_q)}W(k) \cong K$, where 
our Frobenius $\sigma$ on $K$ is a $q$-Frobenius;  
\item[(ii)] $X$ is geometrically connected
\end{enumerate}
in order to prove Theorem \ref{fila} again by Remark \ref{frob} and Lemma \ref{sca}. 

Let $\mathcal M$ be  a convergent $F$-isocrystal on $X/K$ which satisfies the 
hypothesis of Theorem \ref{fila}. Our method of the proof of Theorem \ref{fila} is as follows:
\begin{enumerate}
\item[$1^\circ$] To construct a $\mathrm{Gal}(k(X)^{\mathrm{sep}}/k(X))$-representation $V(\mathcal M)$ 
over $K_\sigma$ which corresponds to the unit-root sub of $M$ at the generic fiber. 
Here $k(X)$ is the function field of $X$ and $k(X)^{\mathrm{sep}}$ is 
a separable closure of $k(X)$. 
\item[$2^\circ$] To show $V(\mathcal M)$ is unramifield at each point of $X$ of codimension $1$. 
Then $V(\mathcal M)$ is a representation of $\pi_1^{\mathrm{et}}(X)$ by Zariski-Nagata purity theorem \cite[X, Theorem 3.4]{SGA2}. 
\item[$3^\circ$] To take a unit-root convergent $F$-isocrystal $\mathcal L$ on $X/K$ by Katz-Crew equivalence 
\cite[Theorem 2.1]{Cr}. 
\item[$4^\circ$] To show $\mathcal L$ is a subobject of $\mathcal M$
\end{enumerate}

\subsection{$(\varphi, \nabla)$-modules over $K[\hspace*{-0.6mm}[t]\hspace*{-0.6mm}]_0$}

Let us assume that the residue field $k$ of $K$ is arbitrary of characteristic $p$ (allowing a non perfect field) 
and $\sigma$ is a $q$-Frobenius on $K$. We also suppose 
there is a $\sigma$-stable Cohen subring $R_0$ of $R$ such that (i) the residue field $R_0$ is $k$, 
(ii) $\mathbb F_q \subset k$ and (iii) $K_\sigma \otimes_{W(\mathbb F_q)}R_0 \cong K$. 

Let us put $K[\hspace*{-0.6mm}[t]\hspace*{-0.6mm}]_0 
= R[\hspace*{-0.6mm}[t]\hspace*{-0.6mm}]\otimes_RK$ 
(resp. $\mathcal E = \widehat{K[\hspace*{-0.6mm}[t]\hspace*{-0.6mm}]_0[1/t]}$ 
to be the $p$-adic completion of $K[\hspace*{-0.6mm}[t]\hspace*{-0.6mm}]_0[1/t]$), 
and $\varphi$ a Frobenius on $K[\hspace*{-0.6mm}[t]\hspace*{-0.6mm}]_0$ 
(resp. the unique extension of $\varphi$ to $\mathcal E$) 
with respect to $\sigma$, that is, 
$\varphi(a) \equiv a^q\, (\mathrm{mod}\, \mathbf m)$ for $a \in R[\hspace*{-0.6mm}[t]\hspace*{-0.6mm}]$ 
and $\varphi|_R = \sigma$. 
We define Gauss norm on $K[\hspace*{-0.6mm}[t]\hspace*{-0.6mm}]_0$ (resp. $\mathcal E$) by 
$$
     \left|\sum_n\, a_nt^n\right|_{\mathrm{Gauss}} = \underset{n}{\mathrm{sup}}\, |a_n|_p, 
$$
where $|a|_p = p^{-\mathrm{ord}_p(a)}$ is the $p$-adic norm. 
Note that $K[\hspace*{-0.6mm}[t]\hspace*{-0.6mm}]_0$ (resp. $\mathcal E$) 
is a principal ideal domain by Weierestrass preparation theorem 
(resp. a completely discrete valuation field). $K[\hspace*{-0.6mm}[t]\hspace*{-0.6mm}]_0$ is complete 
in $(t, p)$-adic topology. 

For either $B = K[\hspace*{-0.6mm}[t]\hspace*{-0.6mm}]_0$ or $\mathcal E$, 
a $(\varphi, \nabla)$-module $(M, \nabla, \Phi)$ over $B$ 
is a free $B$-module 
of finite rank with a $K$-connection $\nabla : M \rightarrow M \otimes_B\Omega_{B/K}^1, \, 
(\Omega_{B/K}^1 = Bdt)$ 
and a Frobenius $\Phi : \varphi^\ast M\, \displaystyle{\mathop{\rightarrow}^{\cong}}\, M$ such that  
$\Phi$ is a horizontal isomorphism with respect to connections. 
Then the category of $(\varphi, \nabla)$-modules over $B$ is Abelian 
and it is independent of the choice of $\varphi$ up to canonical equivalences. 

A generic slope (resp. a special slope) 
of a $(\varphi, \nabla)$-module $M$ over $K[\hspace*{-0.6mm}[t]\hspace*{-0.6mm}]_0$ 
is a slope of the $(\varphi, \nabla)$-module $M \otimes_{K[\hspace*{-0.6mm}[t]\hspace*{-0.6mm}]_0}\mathcal E$ 
over $\mathcal E$ (resp. the $\sigma$-module $M \otimes_{K[\hspace*{-0.6mm}[t]\hspace*{-0.6mm}]_0} K$ 
over $K$, where $K[\hspace*{-0.6mm}[t]\hspace*{-0.6mm}]_0 \rightarrow K$ is defined by $\sum_na_nt^n \mapsto a_0$). 
We also define the Newton polygons at the genetic point and at the special point by the generic slopes 
and the special slopes respectively, 
and the slope filtration on $(\varphi, \nabla)$-module over $\mathcal E$  
and $K[\hspace*{-0.6mm}[t]\hspace*{-0.6mm}]_0$ 
by the same manner if it exists.

\begin{theorem}\label{sfil} 
With the notation above, we have the following. 
\begin{enumerate}
\item \mbox{\rm (\cite[Remark 1.7.8]{Ke0}, \cite[Theorems 2.4]{CT})}
Any $(\varphi, \nabla)$-module $M$ over $\mathcal E$ 
admits a slope filtration $\{ S_\lambda M \}_\lambda$ as $(\varphi, \nabla)$-modules over $\mathcal E$. 
\item \mbox{\rm (\cite[Corollary 2.6.2]{Kat} if $k$ is perfect, 
\cite[Theorems 6.21]{CT} in general)} 
If a $(\varphi, \nabla)$-module $M$ over $K[\hspace*{-0.6mm}[t]\hspace*{-0.6mm}]_0$ 
has constant Newton polygons 
(both generic and special Newton polygons are same), then 
$M$ admits a slope filtration $\{ S_\lambda M \}_\lambda$ as $(\varphi, \nabla)$-modules over 
$K[\hspace*{-0.6mm}[t]\hspace*{-0.6mm}]_0$. 
Moreover, if the residue field $k$ of $K$ is separable closed, then 
$S_0 M \cong (K[\hspace*{-0.6mm}[t]\hspace*{-0.6mm}]_0)^r$ as $(\varphi, \nabla)$-modules 
over $K[\hspace*{-0.6mm}[t]\hspace*{-0.6mm}]_0$. 
\end{enumerate}
\end{theorem}

\vspace*{3mm}

Let $k(\hspace*{-0.5mm}(t)\hspace*{-0.5mm})^{\mathrm{sep}}$ 
(resp. $k^{\mathrm{sep}}$, resp. $k(\hspace*{-0.5mm}(t)\hspace*{-0.5mm})^{\mathrm{ur}}$) 
be a separable closure of $k(\hspace*{-0.5mm}(t)\hspace*{-0.5mm})$, 
(resp. the separable closure of $k$ in $k(\hspace*{-0.5mm}(t)\hspace*{-0.5mm})^{\mathrm{sep}}$, 
resp. the maximal 
unramified extension of $k(\hspace*{-0.5mm}(t)\hspace*{-0.5mm})$ in $k(\hspace*{-0.5mm}(t)\hspace*{-0.5mm})^{\mathrm{sep}}$), 
and $\widehat{\mathcal E}^{\mathrm{nr}}$ (resp. $\widehat{K}^{\mathrm{ur}}$) a $p$-adic completion 
of maximal unramified extension $\mathcal E^{\mathrm{nr}}$ (resp. $K^{\mathrm{ur}}$) 
of $\mathcal E$ (resp. $K$) with residue field $k(\hspace*{-0.5mm}(t)\hspace*{-0.5mm})^{\mathrm{sep}}$ 
(resp. $k^{\mathrm{sep}}$). 
Then $\sigma$ extends on $\widehat{\mathcal E}^{\mathrm{nr}}$ and $\widehat{K}^{\mathrm{ur}}$ uniquely. 
We also use the same notation $\sigma$ on the extensions.
For a $\sigma$-module $M$ over $\mathcal E$, we put 
$$
     V_{\mathcal E}(M) = \mathrm{Ker}(1 -\Phi\otimes\sigma; M \otimes_{\mathcal E} \widehat{\mathcal E}^{\mathrm{nr}}). 
$$
Then $V_{\mathcal E}(M)$ is an $K_\sigma (= \mathcal E_\sigma)$-space 
of dimension less than or equal to $\mathrm{dim}_{\mathcal E}M$ with 
a continuous action of 
$\mathrm{Gal}(k(\hspace*{-0.5mm}(t)\hspace*{-0.5mm})^{\mathrm{sep}}/k(\hspace*{-0.5mm}(t)\hspace*{-0.5mm})) 
= \mathrm{Gal}(\mathcal E^{\mathrm{ur}}/\mathcal E)$ 
by acting $1 \otimes \tau$ for $\tau \in 
\mathrm{Gal}(k(\hspace*{-0.5mm}(t)\hspace*{-0.5mm})^{\mathrm{sep}}/k(\hspace*{-0.5mm}(t)\hspace*{-0.5mm}))$. 
We also define $V_K(N)$ by the same manner for an $F$-space $N$ over $K$. 

\begin{proposition}\label{rep} 
\begin{enumerate}
\item Let $M$ be a $(\varphi, \nabla)$-module over $\mathcal E$ with 
the slope filtration $\{ S_\lambda M\}_\lambda$ such that $S_\lambda M = 0$ for $\lambda < 0$
and $r_0 = \mathrm{dim}_{\mathcal E}\, S_0M$. 
Then $V_{\mathcal E}(M) = V_{\mathcal E}(S_0M)$ is a $K_\sigma$-space of dimension $r_0$ 
with a continuous 
$\mathrm{Gal}(k(\hspace*{-0.5mm}(t)\hspace*{-0.5mm})^{\mathrm{sep}}/k(\hspace*{-0.5mm}(t)\hspace*{-0.5mm}))$-representation. 
\item Let $M$ be a $(\varphi, \nabla)$-module over $K[\hspace*{-0.6mm}[t]\hspace*{-0.6mm}]_0$ with 
the slope filtration $\{ S_\lambda M\}_\lambda$ such that $S_\lambda M = 0$ for $\lambda < 0$. 
Then there are natural isomorphisms
$$
      V_{\mathcal E}(M \otimes_{K[\hspace*{-0.6mm}[t]\hspace*{-0.6mm}]_0} \mathcal E)\, \overset{\cong}{\leftarrow}\,  
      \mathrm{Ker}(1 -\Phi\otimes\sigma; M \otimes_{K[\hspace*{-0.6mm}[t]\hspace*{-0.6mm}]_0} 
      \widehat{K}^{\mathrm{ur}}[\hspace*{-0.6mm}[t]\hspace*{-0.6mm}]_0) \overset{\cong}{\rightarrow}
      V_K(M \otimes_{K[\hspace*{-0.6mm}[t]\hspace*{-0.6mm}]_0} K)
$$
which are compatible with the action of the Galois group 
$\mathrm{Gal}(k(\hspace*{-0.5mm}(t)\hspace*{-0.5mm})^{\mathrm{sep}}/k(\hspace*{-0.5mm}(t)\hspace*{-0.5mm}))$. 
In particular, the representation 
$V_{\mathcal E}(M \otimes_{K[\hspace*{-0.6mm}[t]\hspace*{-0.6mm}]_0} \mathcal E)$ is unramified, that is, 
the Galois group  
acts via the quotient  
$\mathrm{Gal}(k(\hspace*{-0.5mm}(t)\hspace*{-0.5mm})^{\mathrm{ur}}/k(\hspace*{-0.5mm}(t)\hspace*{-0.5mm})) 
\cong \mathrm{Gal}(k^{\mathrm{sep}}/k)$. 
\end{enumerate}
\end{proposition}

\prf{(2) Since the Newton polygons of $M$ both at the generic point and the special point coincide, 
both homomorphisms are bijective by Theorem \ref{sfil}. By the construction 
they are compatible 
with the Galois actions, and hence the representation $V_{\mathcal E}(M \otimes \mathcal E)$ is unramifield. 
}

\subsection{Construction of a functor $V$}

Now we return to the situation of Theorem \ref{fila} with the assumptions at the end of section \ref{filas}. 
We may assume that $X$ is affine and there exists a smooth affine formal 
scheme $\mathrm{Spf}\, A$ 
topologically of finite type over $\mathrm{Spf}\, R$ such that $X = \mathrm{Spec}\, A \otimes_R R/\mathbf m$ 
and $A$ is furnished with a Frobenius endomorphism $\varphi$ compatible to $\sigma$, 
i.e., a continuous ring homomorphism satisfying $\varphi(a)\, \equiv\, a^q\, (\mathrm{mod}\, \mathbf m)$ and 
$\varphi|_R = \sigma$. Such an $A$ exists by \cite[Th\'eor\`eme 6]{El}. 

Let $\mathcal M$ be  a convergent $F$-isocrystal on $X/K$ which satisfies the 
hypothesis of Theorem \ref{fila}. 
Put $M = \Gamma(]X[_{\mathrm{Spf}\, A}, \mathcal M)$. Then $M$  
is a projective $A$-module of finite type which is furnished with an integrable connection 
$\nabla : M \rightarrow M\otimes_A \Omega_{\mathrm{Spf}\, A/\mathrm{Spf}\, R}^1$ 
and a Frobenius $\Phi : \varphi^\ast M\, \displaystyle{\mathop{\rightarrow}^{\cong}}\, M$ which is 
horizontal with respect to integrable connections. 

Let $E$ be a field of fractions of the $p$-adic completion $\widehat{A}_{\mathbf m}$ 
of the localization $A_{\mathbf m}$ at $\mathbf mA$. Then $\widehat{A}_{\mathbf m}/\mathbf m\widehat{A}_{\mathbf m} = k(X)$, 
the field of functions on $X$. 
We define a $K_\sigma$-space ($(\widehat{E}^{\mathrm{ur}})_\varphi = K_\sigma$) by 
$$
    V(\mathcal M) = \mathrm{Ker}(1 -\Phi\otimes\sigma; M \otimes_{A[1/p]}\widehat{E}^{\mathrm{ur}}). 
$$
Then $V(\mathcal M)$ is a $K_\sigma$-space of rank $r_0$ 
with a continuous action of the Galois group 
$$
      \mathrm{Gal}(k(X)^{\mathrm{sep}}/k(X)) = \mathrm{Gal}(E^{\mathrm{ur}}/E). 
$$
Since $\mathcal M$ is a convergent $F$-isocrystal on $X/K$, $V(\mathcal M)$ does not depend on the choice of the Frobenius 
endomorphism $\varphi$ on $\mathrm{Spf}\, A$ and also on the choice of $\mathrm{Spf}\, A$ up to canonical isomorphisms. 

Applying Proposition \ref{rep} at each point of $X$ of codimension $1$, we have 
the following proposition. 

\begin{proposition}\label{unr} 
\begin{enumerate}
\item The Galois group $\mathrm{Gal}(k(X)^{\mathrm{sep}}/k(X))$ 
acts on $V(\mathcal M)$ via the etale fundamental group 
$\pi_1^{\mathrm{et}}(X)$ of $X$. 
\item There is a natural isomorphism $V(\mathcal M\otimes_{\mathcal O_{]X[}}\mathcal N) \cong 
V(\mathcal M) \otimes_{K_\sigma} V(\mathcal N)$ as representations of $\pi_1^{\mathrm{et}}(X)$ 
for convergent $F$-isocrystals $\mathcal M$ and $\mathcal N$ on $X/K$ 
such that both $\mathcal M$ and $\mathcal N$ satisfy the hypothesis $(\ast)$ in Theorem \ref{fila}. 
\end{enumerate}
\end{proposition}

\prf{(1) Let $x$ be a point of $X$ of codimension $1$, and $\mathbf p_x$ the inverse image 
of the prime ideal associated to $x$ by the natural surjection $A \rightarrow A/\mathbf mA$. 
Let $A_x$ be a localization of $A$ at $\mathbf p_x$, $\widehat{A}_x$ the completion of $A_x$ along $\mathbf p_x$, 
$\mathbf p_x\widehat{A}_x = (\mathbf m, t_x)\widehat{A}_x$ for a parameter $t_x \in A_x$, 
and $K_x$ the field of fractions of $A_x/t_xA_x$. 
Such $t_x$ can be taken since $A_x$ is a regular local domain of dimension $2$ 
with the maximal ideal $\mathbf p_x\widehat{A}_x$. 
The ring $\widehat{A}_x$ is integral and we have an isomorphism 
$\widehat{A}_x[1/p] \cong K_x[\hspace*{-0.6mm}[t_x]\hspace*{-0.6mm}]_0$. 
We put $\mathcal E_x$ to be the $p$-adic completion of $K_x[\hspace*{-0.6mm}[t_x]\hspace*{-0.6mm}]_0[1/t_x]$ 
and we regard $E$ as a subfield of $\mathcal E_x$ by the natural injection 
$A \rightarrow \mathcal E_x$. 

Since the natural sequence of isomorphisms 
$$
  V(\mathcal M) \overset{\cong}{\rightarrow} V_E(M \otimes_{A[1/p]} E) 
    \overset{\cong}{\rightarrow}  V_{\mathcal E_x}(M \otimes_{A[1/p]} \mathcal E_x) 
$$
are compatible with Galois actions, $V(\mathcal M)$ is unramified at $x$ 
by Proposition \ref{rep}. Since $x$ is arbitrary of codimension $1$, 
$V(\mathcal M)$ is a representation of $\pi_1^{\mathrm{et}}(X)$ by Zariski-Nagata purity theorem 
\cite[X, Theorem 3.4]{SGA2}. 
}

\vspace*{3mm}

Our functor $V$ is compatible to Katz-Crew equivalence. 

\begin{theorem}\label{crew} Let $X$ be a smooth geometrically connected scheme of finite type over $k$. 
Then the functor $V$ is compatible with the functor $V^0$ of Katz-Crew's equivalence: 
$$
     \begin{array}{ccc}
      F\mbox{\rm -Isoc}(X/K) &\overset{V}{\longrightarrow} &\mbox{\rm Rep}_{K_\sigma}(\pi_1^{\mathrm{et}}(X)) \\
            \bigcup&\underset{V^0}{\nearrow} & \\
           F\mbox{\rm -Isoc}(X/K)^0 
           \end{array}
$$
where $F\mbox{\rm -Isoc}(X/K)^0$ is the full subcategory of 
$F\mbox{\rm -Isoc}(X/K)$ consisting of unit-root $F$-isocrystals and 
$\mbox{\rm Rep}_{K_\sigma}(\pi_1^{\mathrm{et}}(X))$ is the category of 
continuous finite dimensional 
$K_\sigma$-representation of $\pi_1^{\mathrm{et}}(X)$.  
\end{theorem}

\prf{We recall the proof of equivalence of $V^0$ by Crew in \cite[Section 2]{Cr}. 
Let $\mathcal M$ be a unit-root convergent $F$-isocrystal on $X/K$ of rank $r_0$. 
By patching technique, we may assume that $X$ is affine, and hence we follow the notation at the beginning of this 
subsection. 
Let $\mathcal M$ be a unit-root convergent $F$-iscrystal on $X/K$. 
Since $\mathcal M$ is unit-root, 
there is a locally free $\mathcal O_{\mathrm{Spf}\, A}$-module 
of finite rank such that $\mathcal L \otimes_RK = \mathcal M$, 
$\mathcal L$ is stable under the Frobenius $\Phi$ 
and that $\Phi : \varphi^\ast \mathcal L \rightarrow \mathcal L$ 
is an isomorphism of $\mathcal O_{\mathrm{Spf}\, A}$-modules 
(the unit-root condition). 
Such an $\mathcal L$ is called an $F$-lattice $\mathcal M$ and 
it always exists by \cite[Proposition 2.5]{Cr}. 
Then there is a sequence 
$$
     X \leftarrow X_1 \leftarrow X_2 \leftarrow X_3 \leftarrow \cdots
$$
of finite etale Galois coverings of $X$ with structure morphisms $\pi_n : X_n \rightarrow X$ such that 
$\pi^\ast_n(\mathcal L/\mathbf m^n\mathcal L)$ is a trivial $\varphi_n$-module over $A_n$, 
and that the sequence 
$$
     \mathrm{Spf}\, A \leftarrow \mathrm{Spf}\, A_1 \leftarrow \mathrm{Spf}\, A_2 \leftarrow \cdots, 
$$
lifts on $\mathrm{Spf}\, R$ with Frobenius $\varphi_n$. 
Then 
$$
\Gamma_n = \mathrm{Ker}\left(1 - \Phi\otimes \varphi_n; 
\Gamma(\mathrm{Spf}\, A_n, \pi^\ast_n(\mathcal L/\mathbf m^n\mathcal L))\right)
$$ 
is a free $R_\sigma/\mathbf m^nR_\sigma$-module of rank $r_0$ 
with an action of $\mathrm{Gal}(X_n/X)$. 
Here $R_\sigma$ is the integer ring of $K_\sigma$. 
If one puts  
$$
   \begin{array}{l}
     \widehat{B} = \mbox{the $p$-adic completion of}\, \, \cup_n A_n \\
     \widehat\varphi = \mbox{the continuous extension of}\, \, 
     \displaystyle{\mathop{\mathrm{lim}}_{\underset n{\rightarrow}}}\, \varphi_n, 
   \end{array}
$$
then we have 
$$
      V^0(\mathcal M) := K_\sigma \otimes_{R_\sigma} \displaystyle{\mathop{\mathrm{lim}}_{\underset n{\leftarrow}}}\, \Gamma_n 
      \cong K_\sigma \otimes_{R_\sigma} \mathrm{Ker}\left(1 - \Phi \otimes \widehat{\varphi}; 
      \Gamma(\mathrm{Spf}\, A, \mathcal L) \otimes_A\widehat{B} \right) 
      \subset V(\mathcal M)
$$
as $K_\sigma$-spaces with continuous $\pi_1^{\mathrm{et}}(X)$-actions since 
$\widehat{B} \subset \widehat{E}^{\mathrm{ur}}$. 
Comparing the dimension over $K_\sigma$, we have an isomorphism 
$$
V(\mathcal M) \cong V^0(\mathcal M)
$$ 
and Crew's construction coincides with ours. 
}

\subsection{End of the proof of Theorem \ref{fila}} 
Keep the notation as in the previous section. 

\begin{lemma}\label{unitsub} 
Let $\mathcal M$ be a convergent $F$-isocrystal on $X/K$ which satisfies the hypothesis of Theorem \ref{fila}. 
Suppose that $V(\mathcal M)$ contains a trivial representation $V_0 (\cong K_\sigma)$ of rank $1$. Then there exists a 
convergent sub $F$-isocrystal of $\mathcal M$ which is isomorphic to 
the unit-root trivial object $(\mathcal O_{]X[}, d, \varphi)$ of rank $1$. 
\end{lemma}

\prf{We may suppose $X$ is affine and 
the same geometric situation as in the proof of Proposition \ref{unr} holds. 
Let us put $M = \Gamma(]X[_{\mathrm{Spf}\, A}, \mathcal M)$. Then $M$ is a projective $A[1/p]$-module of finite type. 
Our claim is that 
$V_0 \subset M$ under the inclusion $V_0 \subset V(\mathcal M) 
= \mathrm{Ker}(1-\Phi\otimes\varphi; M \otimes_{A[1/p]}\otimes \widehat{E}^{\mathrm{un}}) \subset 
M \otimes_{A[1/p]}\otimes \widehat{E}^{\mathrm{un}}$. 
Since $V_0$ is 
$\mathrm{Gal}(k(X)^{\mathrm{sep}}/k(X))$-invariant, 
$V_0$ is included in $M \otimes_{A[1/p]} E$. 

Let us take a point $x$ of $X$ of codimension $1$, and keep the notation as in the proof of Proposition \ref{unr}. 
Since $V(\mathcal M) =\mathrm{Ker}(1 - \Phi\otimes\sigma; M \otimes_{A[1/p]} \widehat{A}_x[1/p])$ 
by Proposition \ref{rep}, 
$V_0$ is included in $M \otimes_{A[1/p]} \widehat{A}_x[1/p]$. 
Moreover, $M$ is a direct summand of a free $A[1/p]$-module of finite type, the equality 
$$
    M\otimes_{A[1/p]} A_x[1/p] = (M \otimes_{A[1/p]} E) \cap (M \otimes_{A[1/p]} \widehat{A}_x[1/p])
$$
holds in $M \otimes_{A[1/p]} \mathcal E_x$ by Lemma \ref{int} (1). Hence 
$V_0$ is included in $M\otimes_{A[1/p]} A_x[1/p]$. 
Now our claim $V_0 \subset M$ follows from Lemma \ref{int} (2). 
}

\begin{lemma}\label{int} With the notation above, we have
\begin{enumerate}
\item $A_x[1/p] = E \cap \widehat{A}_x[1/p]$. 
\item $A[1/p] = \underset{x}{\cap}\, A_x[1/p]$, where 
$x$ runs through all points $x \in X$ of codimension $1$.
\end{enumerate}
\end{lemma}

\prf{Since $A[1/p]$ is Noetherian normal domain, 
the assertion follows $\mathrm{Spec}\, A_x[1/p] = \mathrm{Spec}\, \widehat{A}_x[1/p]$ 
and $\mathrm{Spec}\, A[1/p] = \coprod_x\mathrm{Spec}\,A_x[1/p]\, \mbox{(a disjoint sum)}$. 
}

\vspace*{3mm}

Now let us complete a proof of Theorem \ref{fila}. 
Let $\mathcal M$ be a convergent $F$-isocrystal on $X/K$ which satisfies the hypothesis of Theorem \ref{fila}. 
Then there is a continuous $K_\sigma$-representation $V(\mathcal M)$ of $\pi_1^{\mathrm{et}}(X)$ of dimension $r$. 
By Katz-Crew's quasi-inverse $D^0$ of $V^0$, 
we have a unit-root convergent $F$-isocrystal $\mathcal L = D^0(V(\mathcal M))$ 
on $X/K$ of rank $r_0$. We have only to prove that there exists a nontrivial homomorphism 
$\mathcal L \rightarrow \mathcal M$. Indeed, by applying this argument to the quotient 
$\mathcal M/\mathrm{Im}(\mathcal L \rightarrow \mathcal M)$ iteratedly, 
one has the maximal quotient $\mathcal N$ of $\mathcal M$ such that 
$\mathcal N$ is of rank $\mathrm{rank}\, \mathcal M - r_0$ 
and it has no slope $0$ sub at the generic point. Then the kernel of $\mathcal M \rightarrow \mathcal N$ 
is a unit-root sub object and it is isomorphic to $\mathcal L$ by Katz-Crew equivalence. 

Applying Proposition \ref{unr} and Lemma \ref{unitsub} to $\mathcal L^\vee \otimes_{\mathcal  O_{]X[}}\mathcal M$, 
we have a unit-root trivial convergent subobject in $\mathcal L^\vee \otimes_{\mathcal  O_{]X[}} \mathcal M$. 
Hence we have a nontrivial homomorphism $\mathcal L \rightarrow \mathcal M$. 
\hspace*{\fill} $\Box$

\end{document}